\documentclass[namedate,webpdf,imanum]{ima-authoring-template}

\usepackage{bm}
\usepackage{subfig}
\usepackage{hyperref}
\usepackage{MnSymbol}
\usepackage{booktabs}
\usepackage{caption}
\usepackage{amsmath,amsthm,amsfonts}
\usepackage{etoolbox}
\usepackage{cprotect}
\usepackage{cleveref}
\usepackage{mathtools}

\crefname{algorithm}{Algorithm}{Algorithms}
\crefname{chapter}{Chapter}{Chapters}
\crefname{equation}{}{Equations}
\crefname{figure}{Fig.}{Figures}
\crefname{table}{Table}{Tables}
\crefname{lemma}{Lemma}{Lemmas}
\crefname{section}{Section}{Sections}
\crefname{subsection}{Section}{Sections}

\makeatletter
\def\myunderbrace#1{\mathop{\vtop{\m@th\ialign{##\crcr
   $\hfil\displaystyle{#1}\hfil$\crcr
   \noalign{\kern3\p@\nointerlineskip}%
   \footnotesize\downbracefill\crcr\noalign{\kern3\p@\nointerlineskip}}}}\limits}
\makeatother

\graphicspath{{../figures/}}

\theoremstyle{thmstyletwo}%
\newtheorem{theorem}{Theorem}
\newtheorem{lemma}[theorem]{Lemma}

\numberwithin{equation}{section}

\begin{document}

\DOI{}
\copyrightyear{}
\vol{}
\pubyear{}
\access{}
\appnotes{}
\copyrightstatement{}
\firstpage{1}


\title[Solving nonlinear ODEs with the ultraspherical spectral method]{Solving nonlinear ODEs with the ultraspherical spectral method}

\author{Ouyuan Qin and Kuan Xu*
\address{\orgdiv{School of Mathematical Sciences}, \orgname{University of Science and Technology of China}, \orgaddress{\street{96 Jinzhai Road}, \postcode{Hefei 230026}, \state{Anhui}, \country{China}}}}

\authormark{O.Qin and K.Xu}

\corresp[*]{Corresponding author: \href{email:kuanxu@ustc.edu.cn}{kuanxu@ustc.edu.cn}}

\received{Date}{0}{Year}
\revised{Date}{0}{Year}
\accepted{Date}{0}{Year}


\abstract{We extend the ultraspherical spectral method to solving nonlinear ODE boundary value problems. We propose to use the inexact Newton-GMRES framework for which an effective preconditioner can be constructed and a fast Jacobian-vector multiplication can be effected, thanks to the structured operators of the ultraspherical spectral method. With a mixed-precision implementation, the inexact Newton-GMRES-ultraspherical framework exhibits extraordinary speed and 
accuracy, as we show by extensive numerical experiments.}


\maketitle

\def\mL{\mathcal{L}}
\def\mB{\mathcal{B}}
\def\bc{\boldsymbol{c}}
\def\md{\mathrm{d}}
\def\bu{\boldsymbol{u}}
\def\mB{\mathcal{B}}
\def\mD{\mathcal{D}}
\def\mF{\mathcal{F}}
\def\mG{\mathcal{G}}
\def\mH{\mathcal{H}}
\def\mI{\mathcal{I}}
\def\mJ{\mathcal{J}}
\def\mK{\mathcal{K}}
\def\mL{\mathcal{L}}
\def\mM{\mathcal{M}}
\def\mN{\mathcal{N}}
\def\mS{\mathcal{S}}
\def\mP{\mathcal{P}}
\def\mQ{\mathcal{Q}}
\def\mO{\mathcal{O}}
\def\mR{\mathcal{R}}
\def\mT{\mathcal{T}}
\def\mW{\mathcal{W}}
\def\ba{\boldsymbol{a}}
\def\bb{\boldsymbol{b}}
\def\bc{\boldsymbol{c}}
\def\be{\boldsymbol{e}_1}
\def\bn{\boldsymbol{n}}
\def\br{\boldsymbol{r}}
\def\bv{\boldsymbol{v}}
\def\lnorm{\left \|}
\def\rnorm{\right \|}
\def\bdlt{\boldsymbol{\delta}}
\def\nl{\left \|}
\def\nr{\right \|}
\def\diag{\mathop{\mathrm{diag}}}
\def\whu{\widehat{u}}
\def\veps{\varepsilon}
\def\udU{\underline {U}}
\def\udA{\underline {A}}
\def\udB{\underline {B}}
\def\eps{\epsilon}
\def\dk{\delta^k}
\def\al{a^{\lambda}}
\def\dl{d^{\lambda}}
\def\ha{\hat{a}}

\section{Introduction} \label{sec:intro}
In this article, we extend the ultraspherical spectral method \citep{olv} to solving the nonlinear ODE boundary value problem
\begin{align*}
\mF(u) = 0, ~\text{s.t.}~~ \mN(u) = 0,
\end{align*}
where $\mF$ is a nonlinear differential operator on $u(x)$. The solution $u(x)$ is a univariate function of the independent variable $x \in [-1, 1]$. By nonlinear, it is meant that $\mF$ cannot be written in the form of \cref{opL} (see \cref{sec:us} below). The functional constraint $\mN$ contains linear or nonlinear boundary conditions or side constraints of other types, such as interior point conditions, global constraints, etc. 

In the very last paragraph of \citet{olv}, the authors briefly discussed the possibility of solving nonlinear differential equations by the ultraspherical spectral method and cautioned the loss of bandedness in the multiplication operators as a threat to the sparsity of the linear system and, therefore, to the exceptional speed of the ultraspherical spectral method. A decade has elapsed since the publication of \citet{olv} and it seems that no progress has been made towards this extension. This paper intends to fill this gap which is long overdue.

What we find in this study is that the loss of bandedness in multiplication operators can be remediated surprisingly easily with the inexact Newton-GMRES framework, a proper preconditioner, and a correct way to carry out the Jacobian-vector multiplication. The proposed framework is extremely easy to implement and can largely restore the lightning speed of the original ultraspherical spectral method without compromising on the accuracy and stability.

We begin by first setting up the inexact Newton-GMRES-ultraspherical (INGU) framework in \cref{sec:ingu}. We discuss the fast application of the truncated Fr\'{e}chet operators in \cref{sec:multiplication} and the preconditioning in \cref{sec:preconditioning}. \cref{sec:acceleration} describes a few possibilities to further speed up the computation, including our mixed precision implementation. Numerical experiments are shown in \cref{sec:experiments} before we close in \cref{sec:conclusion}.

Throughout this article, calligraphy fonts are used for operators or infinite matrices and bold fonts for infinite vectors, whereas functions and the truncations of operators, infinite matrices, and infinite vectors are denoted by normal fonts. We denote the infinite identity operator by $\mI$.

\section{The INGU framework} \label{sec:ingu}
To make our discussion uncluttered, we give a quick review of the very essence of the ultraspherical spectral method in \cref{sec:us}, which we could have dispensed with but at a cost of annoyingly frequent reference to \citet{olv}. We discuss the linearization, truncation, and adaptivity in \cref{sec:linearization} which yield a primitive ultraspherical-based Newton method. The inexact Newton condition enabled by GMRES and the three popular global Newton variants briefly reviewed in \cref{sec:inexact} and \cref{sec:global}, respectively, lead to the prototype INGU method given in \cref{sec:prototype}.

\subsection{Ultraspherical spectral method} \label{sec:us}
The ultraspherical spectral method solves the linear ODE
\begin{align}
\mL u &= f \label{luf}\\
\mbox{s.t. } \mB  u &= c \nonumber
\end{align}
by approximating the solution with a Chebyshev series
\begin{align*}
u(x) = \sum_{k=0}^{\infty}u_k T_k(x),
\end{align*}
where $T_k(x)$ is the Chebyshev polynomials of degree $k$. Here, $\mL$ is an $N$th order linear differential operator
\begin{align}
\mL = a^N(x)\frac{{\md}^N}{{\md}x^N} + \ldots +a^1(x)\frac{{\md}}{{\md}x} + a^0(x), \label{opL}
\end{align}
and $\mB$ contains $N$ linear functionals of boundary conditions. Once the 
coefficients $u_k$ are known, the solution $u(x)$ is identified by the 
coefficient vector $\bu = (u_0, u_1, u_2, \ldots)^{\top}$.

For $\lambda \geq 1$, the $\lambda$th-order differentiation operator has only one nonzero diagonal
\begin{align*}
\mD_{\lambda} = 2^{\lambda - 1}(\lambda - 1)!
\begin{aligned}
  & \quad \myunderbrace{\text{\footnotesize $\lambda$ times}} \\[-7pt]
  & \begin{pmatrix}
       \  0 ~~ \cdots ~~ 0  &\lambda & & & &  \\
      & &\lambda +1&  & & \\
      & & &\lambda +2 & & \\
      & & & &\ddots& \\
  \end{pmatrix}
\end{aligned},
\end{align*}
mapping the Chebyshev $T$ coefficients to the ultraspherical $C^{(\lambda)}$ coefficients\footnote{In \citet{olv}, $\mD_0 = \mD_1$, while in this paper $\mD_1$ maps from Chebyshev $T$ to $C^{(1)}$ and $\mD_0 = \mI$, i.e., the identity operator, for notational consistency.}.

When $a^0(x)= \sum_{k=0}^{\infty}a_j T_j(x)$ is variable, the action of 
$a^0(x)$ on $u(x)$ is represented by a Toeplitz-plus-Hankel-plus-rank-$1$ 
multiplication operator
\begin{align*}
  \mM_{0}[a^0]=\frac{1}{2}\left[\begin{pmatrix}
  2 a_{0} & a_{1} & a_{2} & a_{3} & \cdots \\
  a_{1} & 2 a_{0} & a_{1} & a_{2} & \ddots \\
  a_{2} & a_{1} & 2 a_{0} & a_{1} & \ddots \\
  a_{3} & a_{2} & a_{1} & 2 a_{0} & \ddots \\
  \vdots & \ddots & \ddots & \ddots & \ddots
  \end{pmatrix}+\begin{pmatrix}
  0 & 0 & 0 & 0 & \cdots \\
  a_{1} & a_{2} & a_{3} & a_{4} & \cdots \\
  a_{2} & a_{3} & a_{4} & a_{5} & \udots \\
  a_{3} & a_{4} & a_{5} & a_{6} & \udots \\
  \vdots & \udots & \udots & \udots & \udots
  \end{pmatrix}\right].
\end{align*}
If any of $\al(x) = \sum_{k=0}^{\infty}\al_j C^{(\lambda)}_j(x)$ for $\lambda > 0$ is not constant, the differential operator $\mD_{\lambda}$ should be pre-multiplied by the multiplication operator $\mM_{\lambda}[\al]$, which represents the product of two $C^{(\lambda)}$ series. The most straightforward way to generate $\mM_{\lambda}[\al]$ \citep[Section 6.3.1]{tow} is to express it by the series $\mM_{\lambda}[\al]=\sum_{j=0}^{\infty} a_{j} \mM_{\lambda}\left[C_{j}^{(\lambda)}\right]$, where $\mM_{\lambda}\left[C_{j}^{(\lambda)}\right]$ is obtained by a three-term recurrence relation. Another way to construct these multiplication operators is by an explicit formula given by equation (3.6) in \citet{olv}. Our solution to the nonlinear ODEs relies on a third way to construct $\mM_{\lambda}$, which allows fast applications of a truncation of $\mM_{\lambda}$ to vectors. See \cref{sec:multiplication} below.

When $\mD_{\lambda}$ and $\mM_{\lambda}[\al]$ are employed, each term in 
\cref{opL} maps to a different ultraspherical basis. So the following 
conversion operators $\mS_0$ and $\mS_{\lambda}$ are used to map the 
coefficients in Chebyshev $T$ and $C^{(\lambda)}$ to those in $C^{(1)}$ 
and $C^{(\lambda + 1)}$, respectively

\begin{align*}
  \mS_0 = \begin{pmatrix}
  1&  &-\frac{1}{2}&  & & &\\
   &\frac{1}{2}& &-\frac{1}{2}& & & \\
   & &\frac{1}{2}& &-\frac{1}{2} & & \\
   & & &\ddots& &\ddots& \\
  \end{pmatrix}
  ~~\text{and}~~
  \mS_{\lambda} = \begin{pmatrix}
  1&  &-\frac{\lambda}{\lambda +2}&  & & &\\
   &\frac{\lambda}{\lambda +1}& &-\frac{\lambda}{\lambda +3}& & & \\
   & &\frac{\lambda}{\lambda +2}& &-\frac{\lambda}{\lambda +4} & & \\
   & & &\ddots& &\ddots& \\
  \end{pmatrix}
\end{align*}
for $\lambda \geq 1$. In terms of these operators, the differential operator \cref{opL} can be represented as
\begin{align}
\mL = \mM_N[a^N] \mD_N + \sum_{\lambda = 0}^{N-1}{\mS}_{N-1}\ldots {\mS}_{\lambda}\mM_{\lambda}[\al]\mD_{\lambda}\label{opLinf},
\end{align}
and \cref{luf} becomes
\begin{align}
\mL u = {\mS}_{N-1}\ldots {\mS}_0f, \label{lufInf}
\end{align}
where both the sides are coefficients in $C^{(N)}$.

If $\al(x)$ is analytic or many-times differentiable, it can be approximated to 
machine precision by a Chebyshev approximant of very low degree, resulting in 
the infinite multiplication matrices $\mM_0[a^0]$ and $\mM_{\lambda}[\al]$ 
being banded. Thus, $\al(x)$ is assumed to be a Chebyshev series, and, if the degree of $\al(x)$ is denoted by $\dl$, the bandwidths of $\mM_{\lambda}[\al]$ are 
both $\dl$. Therefore, the $\lambda$th-order operator in the sum in 
\cref{opLinf} has the upper and lower bandwidths $\dl+2N-\lambda$ and 
$\dl-\lambda$, respectively.

To obtain a system of finite dimension, one truncates $\mL$ by premultiplying 
$\mP_{n-N}$ and postmultiplying $\mP_n^{\top}$, where the projection operator 
$\mP_n = ({I}_n , \boldsymbol{0})$.  Incorporating the first $n$ columns of the 
boundary conditions gives us an $n \times n$ square system
\begin{align}
\begin{pmatrix}
  \mB\mP_n^{\top} \\
  \mP_{n-N} \mL \mP_n^{\top}
\end{pmatrix}
\mP_{n}\bu
=
\begin{pmatrix}
c \\
\mP_{n-N}{\mS}_{N-1}\ldots {\mS}_0\mP_{n}^{\top} \mP_{n}\boldsymbol{f}
\end{pmatrix}, \label{lufFin}
\end{align}
where the unknown $\mP_{n}\bu$ and the (unconverted) right-hand side $\mP_{n}\boldsymbol{f}$ are $n$-vectors. Solving \cref{lufFin} gives the Chebyshev coefficients of the $n$-truncation of the solution
\begin{align*}
\widetilde{u}_n(x, t) = \sum_{k=0}^{n-1}u_k T_k(x).
\end{align*}
The matrix on the left-hand side of \cref{lufFin} is almost-banded with the 
upper and lower bandwidths both $N+\displaystyle \max_{\lambda}(\dl-\lambda)$, 
when $n > N+ \displaystyle \max_{\lambda}(\dl-\lambda)$.

The ultraspherical spectral method recapitulated above enjoys three key 
advantages over the collocation-based pseudospectral methods:
\begin{itemize}
\item Since \cref{lufFin} is an almost-banded system, it can be solved in $\mO(n)$ flops, where $n$ is the degrees of freedom of the solution vector.
\item Without preconditioning, the condition number of \cref{lufFin} grows only 
linearly with $n$. When a simple diagonal preconditioner is applied, the 
condition number becomes constant.
\item The forward error of the computed solution can be read directly from the 
right-hand side of the linear system in the course of solving \cref{lufFin} by 
QR factorization. This helps determine the minimal degrees of freedom that is 
required to resolve the solution for a preset accuracy tolerance, therefore 
allowing for adaptivity at virtually no extra cost.\footnote{This is exactly 
how adaptivity is effected in ApproxFun \citep{olv2}.}
\end{itemize}
These advantages are largely retained in solving nonlinear ODEs as we shall see.

\subsection{Linearization, truncation, and adaptivity} \label{sec:linearization}
In $k$th iteration of Newton's method, we solve the linearized problem
\begin{align*}
\mJ[u^k] \dk(x) = -\mF(u^k)
\end{align*}
to obtain the update $\dk(x)$ for the current iterate $u^k(x)$. $\mJ[u^k]$, the 
Fr\'{e}chet derivative of $\mF$ at $u^k(x)$, is a linear differential operator
\begin{align}
\mJ[u^k] = a^N(x)\frac{{\md}^N}{{\md}x^N} + \ldots +a^1(x)\frac{{\md}}{{\md}x} 
+ a^0(x), \label{frechet}
\end{align}
where we slightly abuse the notations by recycling the symbols used in \cref{opL}. Here, we use brackets instead of parentheses in $\mJ[u^k]$ to emphasis that $\mJ[u^k]$ is constructed out of $u^k(x)$ in the sense that the variable coefficients $\al(x)$ depend on $u^k(x)$. It should not be understood as $\mJ$ acting on $u^k(x)$, as it is $\dk(x)$ that $\mJ[u^k]$ acts on. For the definition of the Fr\'{e}chet derivative, see, e.g., \citet[Section 5.3]{atk}. For the calculation of the Fr\'{e}chet derivatives by algorithmic differentiation, see, e.g., \citet{gri,nau,bir1}.

The dependence of $\al(x)$ on $u^k(x)$ is most commonly seen in the form of composition of $u^k(x)$. For example, for the nonlinear operator on the left-hand side of the Bratu equation with a given $\beta$
\begin{align*}
u'' + \beta e^u  = 0, ~~~ s.t. ~~~ u(-1) = u(1) = 0,
\end{align*}
the Fr\'{e}chet derivative
\begin{align*}
\mJ[u] = \displaystyle \frac{{\md}^2}{\md x^2} + \beta e^u
\end{align*}
at a given $u$. Here, $a^0(x) = \beta e^{u(x)}$ is a scaled composition of the 
exponential function and $u(x)$. Multiple approaches are available to calculate 
such a composition. The simplest one is to sample the composition, e.g., 
$e^{u^k(x)}$, on Chebyshev grids of increasing sizes and calculate the 
Chebyshev coefficients by fast cosine transform or FFT until the composition is 
fully resolved. 
With the variable coefficients $\al(x)$ available in its Chebyshev or 
ultraspherical coefficients, we have
\begin{align}
\mJ[u^k] = \mM_N[a^N] \mD_N + \sum_{\lambda = 0}^{N-1}{\mS}_{N-1}\ldots {\mS}_{\lambda}\mM_{\lambda}[\al]\mD_{\lambda}. \label{frechetOPs}
\end{align}

Analogously, linearization gives
\begin{align}
\mN'[u^k]\dk(x) = -\mN(u^k), \label{frechetBCs}
\end{align}
where $\mN'[u^k]$, the Fr\'{e}chet derivative of the boundary condition operator $\mN$ at $u^k(x)$, has dimension $N \times \infty$.

We truncate \cref{frechetOPs} and \cref{frechetBCs} to have
\begin{align}
J_n^k \dk = f^k, \label{truncated}
\end{align}
where $\dk = (\dk_{1}, \dk_{2}, \ldots, \dk_{n})^T$,
\begin{align*}
J_n^k =
    \begin{pmatrix}
      \mN'\mP_n^{\top} \\
      \mP_{n-N} \mJ[u^k] \mP_n^{\top}
    \end{pmatrix}, \text{ and }
f^k = -
    \begin{pmatrix}
      \mN(u^k) \\
      \mP_{n-N}{\mS}_{N-1}\ldots {\mS}_0\mP_{n}^{\top} \mP_{n}\mF(u^k)
    \end{pmatrix}.
\end{align*}
The subscript $n$ in $J_n^k$ is used to indicate the dimension of $J_n^k$. $\mF(u^k)$ is also a composition of $u^k$.

If \cref{truncated} were to be solved by the QR factorization as in \citet{olv}, 
adaptivity would be effected as in the linear case, despite of the loss of 
bandedness. However, since we choose to solve \cref{truncated} using GMRES (see 
below), adaptivity has to be realized in another way. The strategy we follow is 
the one introduced by \citet{aur}. Specifically, the dimension of 
\cref{truncated} is initially determined by the degrees of $\al(x)$ for all 
$\lambda$ and the degree of the residual $\mF(u^k)$. That is, we choose
\begin{align}
n = \max \left(N + \displaystyle \max_{\lambda}(\dl-\lambda), d_{\mF}+1\right), \label{n0}
\end{align}
where $\dl$ is reused to denote the degree of $\al(x)$ in \cref{frechetOPs} and 
$d_{\mF}$ is the degree of $\mF(u^k)$. This explains why \cref{truncated} is 
usually dense or nearly so. Since we are essentially working with functions in 
the Newton-Kantorovich framework \citep{bir2}, it is natural to require the 
emergence of a plateau in $\dk$ before it can be deemed as fully resolved. If 
such a plateau is not seen, we double the dimensions of the system 
to deal with the rapid growth of the high-frequency components due to some of 
the most common nonlinearity, such as $u^2$. Thus, we solve \cref{truncated} 
of the initial size \cref{n0}, check the resolution, augment the system, and repeat until the solution is eventually satisfactorily resolved. Finally, the solution 
is ``chopped'' to trim off the unnecessary trailing coefficients of small 
magnitude. Note that each time when $n$ is doubled the system, the bandedness is restored with the bandwidth half of the dimension of the system. Despite the bandedness, this proportionality still suggests $\mO(n^3)$ operations for solving the system via a direct method.

\subsection{Inexact Newton condition} \label{sec:inexact}

Instead of solving \cref{truncated} exactly, we choose to solve it by enforcing only the inexact Newton condition
\begin{align}
\nl J_n^k \dk - f^k \nr \leq \omega^k \nl f^k \nr, \label{inexact}
\end{align}
where $\omega^k \in [0, 1)$ is the forcing term. The purpose of choosing a 
nonzero $\omega^k$ is to solve \cref{truncated} for $\dk$ to just enough 
precision so that good progress can still be made when far from a solution, but 
also to obtain quadratic convergence when near a solution \citep{eis,eis2,kel1}. This immediately suggests the use of Krylov subspace methods, e.g., GMRES, as these methods produce inexact solutions cheaply. That is, once \cref{inexact} is satisfied, 
the iteration of GMRES terminates and returns an inexact solution $\dk$. The 
adoption of the inexact Newton condition and the use of GMRES are justified in 
\cref{sec:preconditioning,sec:multiplication,sec:acceleration} where we see how 
the full potential of the ultraspherical spectral method can be unleashed in 
the current context, forming an extremely efficient framework for solving 
nonlinear ODE boundary value problems.

\subsection{Global Newton methods} \label{sec:global}
A naive implementation of Newton's method based on \cref{truncated}, sometimes 
regarded as the local Newton method, has limited chance of convergence unless 
the initial iterate is close enough to the solution. Thus, a practical global 
Newton method must be applied. The global Newton method for solving nonlinear 
systems has many variants, and they differ mainly by the strategy for 
determining the search direction and the step length, how the linear systems 
are solved, and whether the derivative is dispensed with. We use three global 
Newton methods as the vehicle for implementing the INGU framework -- the trust 
region method \citep[Section 11.2]{noc} enabled by the dogleg approximation 
(TR-dogleg) \citep{pow}, the line search with Armijo backtracking 
(LS-backtracking) \citep[Section 8]{kel1}, and the trust region method in the affine 
contravariant framework (TR-contravariant) \citep[Section 3.2]{deu}.

The TR-dogleg method is arguably the most reliable and widely accepted 
algorithm, being the default implementation in many mainstream computing 
platforms, standard libraries, or public domain codes. However, there is a catch 
-- the TR-dogleg requires the knowledge of the transpose of the Jacobian or the 
access to it via its products with vectors. The line search method is usually 
less robust but relatively easy to implement and this is true particularly for 
its backtracking-based implementation. The TR-contravariant method assumes that 
$\mF(u)$ satisfies an affine contravariant Lipschitz condition based on which 
the minimization of the residual is modeled as a constrained quadratic 
optimization problem. Though it is not as well received as the TR-dogleg and 
LS-backtracking methods, TR-contravariant performs equally well on average in 
our experiments.

\subsection{A prototype framework}\label{sec:prototype}
We now have the key ingredients for setting up the prototype INGU framework, 
which is summarized in \cref{algo:proto}. For convenience, we denote by $\mG$ 
the operator formed by concatenating $\mN$ and $\mF$ vertically as 
$\left[\mN;\mF\right]$, where \textsc{Matlab} syntax is used. This way, when 
$\mG$ is applied to a solution $u^k$, it returns a column vector containing the 
residual of both the nonlinear equation and the boundary conditions.

\begin{algorithm}[!t]
\caption{INGU prototype}  \label{algo:proto}
\begin{unlist}
\item[Inputs: ] nonlinear operator $\mG$, Fr\'{e}chet operator $\mJ$, initial iterate $u^0$, relative tolerance $\eta_r$ for the residual.
\item[Output: ] approximate solution $u^k$ satisfies $\left\|\mG(u^k)\right\| \leq \eta_r \left\|\mG(u^0)\right\| + \eta_r$.
\end{unlist}
\hrule
\begin{algorithmic}[1]
    \State Set $k=0$, construct $u^0$, $\eta = \eta_r \left\|\mG(u^0)\right\| + \eta_r$.
    \While{$\left\|\mG(u^k)\right\| > \eta$}  \Comment{outer iteration} \label{outer}
    \While{$\dk$ is not resolved}  \Comment{intermediate iteration}
    \State Solve \cref{truncated} inexactly by GMRES. \Comment{inner iteration by mixed precision}
    \State Double the size of the system.
    \EndWhile
    \State Call \textsc{Postprocess}.
    \State $k = k+1$
    \EndWhile
    \State \Return $u^k$
\end{algorithmic}
\end{algorithm}

There are three iterations in this framework. As the outer loop (lines 2-9), 
the Newton iteration generates the sequence of the updates $\dk$ and the 
approximate solutions $u^k$. The initial iterate is usually chosen to be the 
polynomial of lowest degree that satisfies the boundary conditions. The outer 
iteration terminates when the residual $\mG(u^k)$ is smaller than the preset 
tolerance. This tolerance usually includes terms for both the relative and 
absolute residuals.

The intermediate loop (lines 3-6) ensures that the update $\dk$ has an adequate 
resolution by keeping doubling the size of the linear system until a plateau is 
formed. A zero vector is usually used as the initial iterate of the 
intermediate loop. The solution from the previous iteration is elongated with 
zeros before fed into the next iteration. In most of our experiments, $n$ is 
doubled only a couple of times before the intermediate loop terminates. See 
\cref{sec:experiments} for more detail.

The Krylov subspace iteration of GMRES is deemed as the inner loop (line 4) and 
is terminated when \cref{inexact} is satisfied. The forcing term $\omega^k$ is 
determined by the function \textsc{Postprocess} from the previous outer 
iteration. The function \textsc{Postprocess} usually takes in the information 
of the current outer iteration, such as the residual and the Jacobian, and the 
information inherited from the previous outer iteration, such as the 
contraction factor and the contravariant Kantorovich quantity for the 
TR-contravariant method. The outputs of \textsc{Postprocess} usually include a 
updated solution $u^{k+1}$ obtained by adding to $u^k$ a post-processed Newton 
step $\widetilde{\delta}^k$, a new forcing term $\omega^{k+1}$, and the size of 
the trust region for the next outer iteration. See 
\cref{algo:postTRDogleg,algo:postLSBacktracking,algo:postTRContravariant} for 
details of the function \textsc{Postprocess} for each global Newton method.

\section{Fast matrix-vector multiplication} \label{sec:multiplication}
The efficiency of the GMRES method hinges on if fast matrix-vector multiplication is available. Specifically, we wish to be able to apply $J_n^k$ to a given $n$-vector speedily, despite the loss of the bandedness in $J_n^k$.

For the moment, let us ignore the top $N$ rows of $J_n^k$, i.e., the boundary conditions, and examine $\mP_{n-N} \mJ \mP_n^{\top}$ in a termwise manner. We first look at the zeroth-order term which is a truncation of ${\mS}_{N-1}\ldots {\mS}_0\mM_0[a^0]$. Denote by $\mT[a^0]$, $\mH[a^0]$, and $\mR[a^0]$ the Toeplitz, the Hankel, and the rank-$1$ parts of $2\mM_0[a^0]$, respectively, i.e., $\mM_0[a^0]=\left(\mT[a^0]+\mH[a^0]+\mR[a^0]\right)/2$, where
\begin{align*}
\mT[a^0] = \begin{pmatrix}
  2 a_0 & a_1 & a_2 & a_3 & \cdots \\
  a_1 & 2 a_0 & a_1 & a_2 & \ddots \\
  a_2 & a_1 & 2 a_0 & a_1 & \ddots \\
  a_3 & a_2 & a_1 & 2 a_0 & \ddots \\
  \vdots & \ddots & \ddots & \ddots & \ddots
  \end{pmatrix},~~~
\mH[a^0]=\begin{pmatrix}
  a_0 & a_1 & a_2 & a_3 & \cdots \\
  a_1 & a_2 & a_3 & a_4 & \udots \\
  a_2 & a_3 & a_4 & a_5 & \udots \\
  a_3 & a_4 & a_5 & a_6 & \udots \\
  \vdots & \udots & \udots & \udots & \udots
  \end{pmatrix},
\end{align*}
and $\mR[a^0] = -\be \ba^0$ with $\be$ the first column of $\mI$ and $\ba^0$ is 
the infinite vector obtained by prolonging $a^0=(a^0_0, a^0_1, \dots, 
a^0_{d^0})$ with zeros. For an $n$-vector $v$, calculating each of $s_T = 
\mP_n\mT[a^0]\mP_n^{\top} v$ and $s_H = \mP_n\mH[a^0]\mP_n^{\top} v$ costs two 
FFTs and one inverse FFT, all of length $2n-1$ \citep[P4.8.6]{gol}. Calculating 
$s_R = 
\mP_n\mR[a^0]\mP_n^{\top} v$, $s = s_T+s_H+s_R$, and $\mP_{n-N}{\mS}_{N-1}\ldots 
{\mS}_0\mP_n^{\top} s$ can be done in $\mO(n)$ flops. Hence, the total cost of 
applying $\mP_{n-N} {\mS}_{N-1}\ldots {\mS}_0\mM_0[a^0] \mP_n^{\top}$ is 
dominated by the six FFTs.

The first-order term $\mP_{n-N} {\mS}_{N-1}\ldots {\mS}_1 \mM_1[a^1] \mD_1 \mP_n^{\top}$ can also be quickly applied, as pointed out in \citet[Remark 3]{olv}. To see this, we observe that $\mM_1[a^1]$ is constructed using $a^1(x)$'s $C^{(1)}$ coefficients $a^1=(a^1_0, a^1_1, \dots, a^1_{d^1})$ and it also acts on and maps to coefficient vectors in $C^{(1)}$. This suggests another way to express it
\begin{align*}
\mM_1[a^1] = {\mS}_0 \mM_0[\mP_{d^1+1}{\mS}^{-1}_0\mP_{d^1+1}^{\top} a^1]{\mS}^{-1}_0, 
\end{align*}
where $\mM_0[\mP_{d^1+1}{\mS}^{-1}_0\mP_{d^1+1}^{\top} a^1]$ means to first 
convert the $C^{(1)}$ coefficients of $a^1(x)$ to the Chebyshev ones and then 
construct the $\mM_0$ multiplication matrix with the Chebyshev coefficients of 
$a^1(x)$. In practice, $a^1(x)$ is often available by its Chebyshev 
coefficients, either as the variable coefficient of a linear term in $\mF(u)$ 
or a composition of $u^k(x)$, or a combination of both. Hence, we assume 
$\ha^1=\mP_{d^1+1}{\mS}^{-1}_0\mP_{d^1+1}^{\top} 
a^1=(\ha_0,\ha_1,\dots,\ha_{d^1})$ is available from now on. The Sylvester map 
or the displacement of $\mM_0[\ha^1]$ is
\begin{align*}
\nabla_{\mS_0}\left(\mM_0[\ha^1]\right) &= \mS_0 \mM_0[\ha^1] - \mM_0[\ha^1] \mS_0 \\
&= \frac{1}{2} \left( \nabla_{\mS_0}\left(\mT[\ha^1]\right) + \nabla_{\mS_0}\left(\mH[\ha^1]\right) + \nabla_{\mS_0}\left(\mR[\ha^1]\right) \right),
\end{align*}
where $\mT[\ha^1]$, $\mH[\ha^1]$, and $\mR[\ha^1]$ are the Toeplitz, the Hankel, and the rank-$1$ parts of $\mM_0[\ha^1]$, respectively. Some algebraic work then gives
\begin{align}
\nabla_{\mS_0}\left(\mT[\ha^1]\right) = \bc_{t} \br_{th}, ~~~\nabla_{\mS_0}\left(\mH[\ha^1]\right) = \hat{\mH}[\ha^1] + \bc_{h} \br_{th}, ~~~\nabla_{\mS_0}\left(\mR[\ha^1]\right) = \bc_{r} \br_{r}, \label{thr}
\end{align}
where
\begin{align*}
& \bc_{t} = \frac{1}{2}
\begin{pmatrix}
2 & 0 & 0 & 0 &\cdots \\
\ha_0+\ha_2 & \ha_1+\ha_3 & \ha_2+\ha_4 & \ha_3+\ha_5 & \cdots \\
\ha_1 & \ha_2 & \ha_3 & \ha_4 & \cdots
\end{pmatrix}^{\top},
\br_{th} = \frac{1}{2}
\begin{pmatrix}
\ha_0 & \ha_1 & \ha_2 & \ha_3 & \cdots \\
-2 & 0 & 0 & 0 & \cdots \\
0 & -2 & 0 & 0 & \cdots \\
\end{pmatrix}, \\
& \bc_{h} = \frac{1}{2}
\begin{pmatrix}
    2 & 0 & 0 & 0 & \cdots \\
    \ha_{\lvert -2 \rvert}+\ha_0 & \ha_{\lvert -1 \rvert}+\ha_1 & \ha_0+\ha_2 & \ha_1+\ha_3 &\cdots \\
    \ha_{\lvert -1 \rvert} & \ha_0 & \ha_1 & \ha_2 & \cdots
\end{pmatrix}^{\top}, \\
& \bc_{r} = -\begin{pmatrix}
    1 & 0 & 0 & 0 & \cdots
  \end{pmatrix}^{\top}, \quad
  \br_{r} = \frac{1}{2}\begin{pmatrix}
    0 & \ha_1 & \ha_0 + \ha_2 & \ha_1 + \ha_3 & \cdots \\
  \end{pmatrix}, \\
&  \hat{\mH}[\ha^1] = \frac{1}{2}\begin{pmatrix}
    \ha_{\lvert -2 \rvert} - \ha_2 & \ha_{\lvert -1 \rvert} - \ha_3 & \ha_0 - \ha_4 & \ha_1 - \ha_5 & \cdots \\
    \ha_{\lvert -1 \rvert} - \ha_3 & \ha_0 - \ha_4 & \ha_1 - \ha_5 & \ha_2 - \ha_6 & \udots \\
    \ha_0 - \ha_4 & \ha_1 - \ha_5 & \ha_2 - \ha_6 & \ha_3 - \ha_7 & \udots \\
    \ha_1 - \ha_5 & \ha_2 - \ha_6 & \ha_3 - \ha_7 & \ha_4 - \ha_8 & \udots \\
    \vdots & \udots & \udots & \udots & \udots \\
  \end{pmatrix}.
\end{align*}
Here, $\hat{\mH}[\ha^1]$ is, again, a Hankel matrix. Absolute values are used in some of the subscripts in $\bc_h$ and $\hat{\mH}[\ha^1]$ to reveal the pattern. Equation \cref{thr} shows that the displacement $\nabla_{\mS_0}\left(\mM_0[\ha^1]\right)$ is a Hankel-plus-low-rank matrix. Expressing $\mM_1[\ha^1]$ in terms of $\nabla_{\mS_0}\left(\mM_0[\ha^1]\right)$, we have
\begin{align}
\mM_1[\ha^1] &= \nabla_{\mS_0}\left(\mM_0[\ha^1]\right){\mS}^{-1}_0 + \mM_0[\ha^1] \nonumber \\ &= \frac{1}{2} \left(\left( \hat{\mH}[\ha^1] + (\bc_{t} + \bc_{h}) \br_{th} + \bc_{r} \br_{r}\right){\mS}^{-1}_0 + \mT[\ha^1] + \mH[\ha^1] + \mR[\ha^1] \right). \label{M1ahat1}
\end{align}
Using \textsc{Matlab}'s notation \texttt{triu} to denote the upper triangular 
part of a matrix, we have
\begin{align}
\mS_{0}^{-1} = \texttt{triu} \left(  \bc^{0} \br \right), ~~~
\bc^{0} = \begin{pmatrix}
    1 & 0 & 2 & 0 & 2 & \cdots \\
    0 & 2 & 0 & 2 & 0 & \cdots \\
\end{pmatrix}^{\top}, ~~~
\br = \begin{pmatrix}
    1 & 0 & 1 & 0 & 1 & \cdots \\
    0 & 1 & 0 & 1 & 0 & \cdots \\
\end{pmatrix}. \label{S0inv}
\end{align}
With \cref{S0inv}, the first term in the outermost parentheses of  \cref{M1ahat1} can be simplified as
\begin{align*}
\left( \hat{\mH}[\ha^1] + (\bc_{t} + \bc_{h}) \br_{th} + \bc_{r} \br_{r}\right){\mS}^{-1}_0
= -\begin{pmatrix}
    \ha_2 & \ha_3 & \ha_4 & \ha_5 & \cdots \\
    \ha_3 & \ha_4 & \ha_5 & \ha_6 & \udots \\
    \ha_4 & \ha_5 & \ha_6 & \ha_7 & \udots \\
    \ha_5 & \ha_6 & \ha_7 & \ha_8 & \udots \\
    \vdots & \udots & \udots & \udots & \udots \\
  \end{pmatrix} - \mH[\ha^1] - \mR[\ha^1].
\end{align*}
The last equation and \cref{M1ahat1} show that $\mM_1[\ha^1]$ is a Toeplitz-plus-Hankel operator given by
\begin{align*}
\mM_1[\ha^1] = \frac{1}{2}
\begin{pmatrix}
2\ha_0 & \ha_1 & \ha_2 & \ha_3 & \cdots \\
\ha_1 & 2\ha_0 & \ha_1 & \ha_2 & \ddots \\
\ha_2 & \ha_1 & 2\ha_0 & \ha_1 & \ddots \\
\ha_3 & \ha_2 & \ha_1 & 2\ha_0 & \ddots \\
\vdots & \ddots & \ddots & \ddots & \ddots \\
\end{pmatrix} -
\frac{1}{2}
\begin{pmatrix}
\ha_2 & \ha_3 & \ha_4 & \ha_5 & \cdots \\
\ha_3 & \ha_4 & \ha_5 & \ha_6 & \udots \\
\ha_4 & \ha_5 & \ha_6 & \ha_7 & \udots \\
\ha_5 & \ha_6 & \ha_7 & \ha_8 & \udots \\
\vdots & \udots & \udots & \udots & \udots \\
\end{pmatrix}. 
\end{align*}
Hence, when the Chebyshev coefficients of $a^1(x)$ are available, the cost of multiplying $\mP_{n-N} {\mS}_{N-1}\ldots {\mS}_1 \mM_1[a^1] \mD_1 \mP_n^{\top}$ with a vector is again six FFTs plus $\mO(n)$ flops.

We wish the higher order terms in $\mM[\al]$ are structured alike so that fast 
multiplications can be effected similarly. Unfortunately, whether general 
higher order multiplication operators bear a Toeplitz plus Hankel form, or 
something akin to to allow a fast application is not known. Fortunately, we can
always circumvent this using the idea above for $\mM[a^1]$. Keeping in mind 
that $\mM_{\lambda}[\al]$ is constructed by using $\al(x)$'s $C^{(\lambda)}$
coefficients $\al = (\al_0, \al_1, \dots, \al_{d^{\lambda}})$ and that it maps 
between $C^{(\lambda)}$ coefficients, we re-express it as
\begin{align}
\mM_{\lambda}[\al] = {\mS}_{\lambda-1}\dots {\mS}_1 \mM_1[\ha^{\lambda}]{\mS}^{-1}_1 \dots {\mS}^{-1}_{\lambda-1}, \label{Ml2}
\end{align}
where $\ha^{\lambda} = \mP_{\dl+1}{\mS}^{-1}_0 \dots {\mS}^{-1}_{\lambda-1} \mP_{\dl+1}^{\top} \al$ are the Chebyshev coefficients of $\al(x)$ and $\mS_{\lambda}^{-1}$ can be expressed explicitly as
\begin{align*}
\mS_{\lambda}^{-1} = \texttt{triu} \left(  \bc^{\lambda} \br \right),~~ \bc^{\lambda} = \frac{1}{\lambda}\begin{pmatrix}
\lambda & 0 & \lambda+2 & 0 & \lambda+4 & \cdots \\
0 & \lambda+1 & 0 & \lambda+3 & 0 & \cdots \\
\end{pmatrix}^{\top}, \text{ for } \lambda \geq 1.
\end{align*}
Substituting \cref{Ml2} into \cref{frechetOPs} yields
\begin{align}
\mJ[u^k] = {\mS}_{N-1}\dots {\mS}_1 \left(\sum_{\lambda = 1}^N \mM_1[\ha^{\lambda}]{\mS}^{-1}_1 \dots {\mS}^{-1}_{\lambda-1} \mD_{\lambda} + \mS_0\mM_0[a^0]\right), \label{fastfrechet}
\end{align}
which shows that the multiplication part of each term in $\mJ[u^k]$ can be done 
via $\mM_1$ or $\mM_0$. Hence, applying $\mP_{n-N}\mJ[u^k]\mP_n^{\top}$ costs 
$6(N+1)$ FFTs. See \cref{algo:fast} for details and the stepwise 
costs. 

\begin{algorithm}[!t]
\caption{Fast Jacobian-vector multiplication}
\label{algo:fast}
\begin{unlist}
\item[Inputs: ] A Fr\'{e}chet operator $\mJ[u^k]$ in the form of \cref{frechetOPs}, the order $N$ of $\mJ[u^k]$, and an $n$-vector $v$.
\item[Output: ] An $n$-vector $v_s = \mP_{n-N}\mJ[u^k]\mP_n^{\top} v$.
\end{unlist}
\hrule
\begin{algorithmic}[1]
\State Calculate the Chebyshev coefficients of $\ha^{\lambda}$ for all 
$\lambda$. \Comment{$\mO(n \log_2 n)$}
\State Calculate $v_s = \mP_n \mS_0\mM_0[a^0] \mP_n^{\top} v$. \Comment{$\mO(n \log_2 n)$}
\For {$\lambda = 1$ to $N$}
\State Calculate $v_s = v_s + \mP_n \mM_1[\ha^{\lambda}]{\mS}^{-1}_1 \dots {\mS}^{-1}_{\lambda-1} \mD_{\lambda} \mP_n^{\top} v$.
\Comment{$\mO(n \log_2 n)$}
\EndFor
\State Calculate and return $\mP_{n-N}{\mS}_{N-1}\dots {\mS}_1 \mP_n^{\top} v_s$.
\Comment{$\mO(n)$}
\end{algorithmic}
\end{algorithm}

Standard FFT libraries, like FFTW \citep{fri}, allow the users to pre-plan an 
optimized FFT of a given size and apply the plan repeatedly to vectors of 
the same size. This can help accelerate each of the intermediate iteration 
for the dimension of the system is unchanged throughout.

Noting that applying the boundary rows $\mN'\mP_n^{\top}$ to $v$ can be done in 
$\mO(Nn)$ flops, we conclude that the multiplication of $J_n^k$ and $v$ costs 
only $\mO(Nn\log_2 n)$ flops. This justifies the use of GMRES.

We also remark that computing $f^k$ costs at most $\mO(n\log_2n)$ flops, since 
$\mF(u^k)$ is a composition of $u^k$ and $\mN(u^k)$ are $N$ functionals.

\section{Preconditioner} \label{sec:preconditioning}
Our GMRES-based approach is also justified by a simple but effective 
preconditioner. The fact that $\mJ[u^k]$ is dense motivates us to use an 
almost-banded preconditioner  --- if a diagonal scaling or Jacobi-type 
preconditioner works perfectly for an almost-banded system, as suggested in 
\citet{olv}, why not use an almost-banded one to precondition the dense system 
which can be deemed as obtained from the same almost-banded system by extending 
the bandwidth to the full dimension of the system? Hence, we propose the use of 
a right preconditioner for the $k$th outer iteration
\begin{align*}
\mW^k =
    \begin{pmatrix}
      \mN' \\
      \widetilde{\mJ}[u^k]
    \end{pmatrix},
\end{align*}
where
\begin{align}
\widetilde{\mJ}[u^k] = \mM_N[\widetilde{a}^N] \mD_N + \sum_{\lambda = 0}^{N-1}{\mS}_{N-1}\ldots {\mS}_{\lambda}\mM_{\lambda}[\widetilde{a}^{\lambda}]\mD_{\lambda}. \label{frechetBanded}
\end{align}
Here, the multiplication operators are constructed from the first 
$m^{\lambda}+1$ leading coefficients of $\al$, that is,
\begin{align*}
\widetilde{a}^{\lambda} = \left(\al_0, \al_1, \ldots, \al_{m^{\lambda}}\right),
\end{align*}
where $m^{\lambda} = p + \lambda \ll n$ and the value of integer $p$ is to be 
determined. The operators $\mM_N[\widetilde{a}^N] \mD_N$ and ${\mS}_{N-1}\ldots 
{\mS}_{\lambda}\mM_{\lambda}[\widetilde{a}^{\lambda}]\mD_{\lambda}$ in 
\cref{frechetBanded} and, therefore, $\widetilde{\mJ}[u^k]$ have the upper and 
lower bandwidths $p+2N$ and $p$ respectively, due to the argument given below 
\cref{lufInf}. Hence, $\widetilde{\mJ}[u^k]$ is a banded approximation of 
$\mJ[u^k]$ in the sense that every component of $\widetilde{a}^{\lambda}$ has 
its contribution in all the nonzero diagonal entries of $\widetilde{\mJ}[u^k]$. 
Note that an entry in the band of $\widetilde{\mJ}[u^k]$ differs from the entry 
in the same position in $\mJ[u^k]$ as the latter has contribution from every 
component of $\al$, not just the first $m^{\lambda}+1$ ones.

Instead of \cref{frechetBanded}, we find it easiest to follow \cref{fastfrechet} to construct the banded part of $\mW^k$. That is,
\begin{align*}
\check{\mJ}[u^k] = {\mS}_{N-1}\dots {\mS}_1 \left(\sum_{\lambda = 1}^N \mM_1[\check{a}^{\lambda}]{\mS}^{-1}_1 \dots {\mS}^{-1}_{\lambda-1} \mD_{\lambda} + \mS_0\mM_0[\check{a}^0]\right),
\end{align*}
where $\check{a}^{\lambda}=\mP_{m^{\lambda}+1}{\mS}^{-1}_0 \dots {\mS}^{-1}_{\lambda-1}\mP_{m^{\lambda}+1}^{\top} \widetilde{a}^{\lambda}$ are the Chebyshev coefficients of $\widetilde{a}^{\lambda}(x) = \sum_{k=0}^{\dl}\widetilde{a}^{\lambda}_j C^{(\lambda)}_j(x)$. Because of the bandedness, the equivalence of $\check{\mJ}[u^k]$ and $\widetilde{\mJ}[u^k]$ can be guaranteed by exact truncations. Instead of solving \cref{truncated}, we solve
\begin{align}
J_n^k \left(W_n^k\right)^{-1} \theta^k = f^k, \label{preconditioned}
\end{align}
where $W_n^k = \mP_n \mW^k \mP_n^{\top}$, and $\dk$ is finally recovered by solving $W_n^k \dk = \theta^k$. Note that $W_n^k$ stays unchanged within each inner loop. Thus, it suffices to compute the QR factorization of $W_n^k$ only once for each call of GMRES. Since the construction, the application, and the inversion of $\check{\mJ}[u^k]$ all cost $\mO(p^2n)$ flops, we choose
\begin{align}
p = \left\lfloor\sqrt{\log_2 n}\right\rfloor \label{p}
\end{align}
to match up to the cost of the FFT-based Jacobian-vector multiplication and the function composition of $u^k$, resulting in an asymptotic complexity of $\mO(n \log_2 n)$. For a nonlinear ODE that is not singularly perturbed, as we shall see in 
\cref{sec:experiments}, $p$ usually has a value below $10$. Our experiments 
suggest that a fixed $p$ of small integral value often works equally well. But 
\cref{p} offers a weak dependence on $n$ and this adaptivity may play a bigger 
role when we migrate to nonlinear problems in higher spatial dimensions.

Assuming $a^N(x)=1$ in \cref{frechet} and employing virtually the same 
technique used in the proof of Lemma 4.3 in \citet{olv}, we can readily show 
that the right-preconditioned system \cref{preconditioned} is a compact 
perturbation of the identity operator in the Banach space $\ell_{K}^2$. For the 
definition of the norm of $\ell_K^2$, see \citet[Definition 4.2]{olv}.

\begin{lemma} \label{thm:compact}
Assume the boundary operator $\mN':\ell_D^2 \rightarrow \mathbb{C}^N$ is 
bounded and $a^N(x)=1$. Let $\mW^k:\ell_{K+1}^2 \rightarrow \ell_K^2$ for some 
$K \in \{ D-1, D, \ldots\}$, where $D$ is the smallest integer such that 
$\mN':\ell^2_D \rightarrow \mathbb{C}^N$ is bounded. Then
\begin{align*}
\begin{pmatrix}
\mN' \\
\mJ[u^k]
\end{pmatrix} \left(\mW^k\right)^{-1}= \mI+\mK,
\end{align*}
where $\mK:\ell_{K}^2 \rightarrow \ell_{K}^2$ is a compact operator for $K = D-1,D,\ldots$.
\end{lemma}

The well-conditionedness implied by \cref{thm:compact} follows from Lemma 4.4 of \citet{olv} and is confirmed in  \cref{sec:eigen} by extensive numerical experiments.

There are two remarks to be made regarding the proposed preconditioner. First, 
$\mW$ can be interpreted as a \emph{coarse-grid} preconditioner or 
\emph{low-order discretization} preconditioner which captures the low-frequency 
components of the problem, leaving the high frequencies to be treated by the 
Krylov subspace iteration. What is also noteworthy is that it works in the 
frequency/coefficient space directly --- there is no need to do interpolation 
and transfer back and forth between the physical/value and the 
frequency/coefficient spaces. 

Second, the diagonal preconditioner $\mR$ given in \citet[Section 4.1]{olv} 
continues to work, despite the loss of bandedness in \cref{truncated}. On the 
one hand, it is apparent that the diagonal preconditioner costs less to apply. 
On the other, experiments show that the proposed preconditioner has a better 
chance to make the eigenvalues of $J_n^k$ cluster\footnote{Of course, 
eigenvalues clustering needs not indicate fast convergence \citep{gre1,gre2}.}. 
Therefore, it is difficult to say which preconditioner is more effective. 
Whether the proposed preconditioner (significantly) outperforms the diagonal 
one also depends on other factors, including but not limited to how frequently 
GMRES is restarted, the forcing term $\omega^k$, etc. Extensive numerical 
tests show that the INGU method is faster with the new preconditioner and very 
much so especially when the number of iterations allowed before GMRES restarts 
is not very large.



\section{Further acceleration} \label{sec:acceleration}
The fast application of the Jacobians and the preconditioner are decisive in 
making our INGU method fast. In this section, we discuss other opportunities 
that may potentially allow the computation to be further accelerated.

\subsection{Mixed precision} \label{sec:mixed}
After a rapid development in the last couple of decades, mixed precision algorithms have earned a proven track record in accelerating iterative methods and made inroads into the tool set of our day-to-day computation \citep{abd, hig}.

Tisseur analyzes the limiting accuracy and limiting residual of Newton's method 
in floating point arithmetic in a multiple-precision setting \citep{tis}. A 
recent work by \citet{kel3} investigates the use of reduced precision 
arithmetic to solve the linearized equation for the Newton update by a direct 
linear solver. These works are reviewed and summarized in \citet[Section 5]{hig}. See 
also Algorithm 5.1 therein. The main idea is to compute the residual in a 
relatively high precision $p_h$ and the Newton update in a relatively low 
precision $p_l$, while maintain the approximate solution $u^k$ in a working 
precision $p_w$ with $p_h \leq p_w \leq p_l$. In \citet{kel3}, $p_h = p_w$ are 
chosen to be double precision and $p_l$ single or half precisions.

Our mixed-precision implementation for the INGU method follows the same 
strategy. To be specific, we compute and store the residual and the solution in 
double precision, whereas the GMRES solve is done in single precision. Similar 
to what is reported in \citet{kel3}, our numerical experiments show that the 
results from a reduced-precision implementation of GMRES does not distinguish 
from those done by a fixed-precision computation with double precision 
throughout. We observe an average gain of $20\%$ to $30\%$ in speed.

There are a few more opportunities for further GMRES speed boost by mixed-precision arithmetic. First, it has been shown in \citet{van,sim,gir} that matrix-vector products can be done in increasingly reduced precision without degrading the overall accuracy of the entire computation. In addition, the orthonormalization step in GMRES can also be performed in a reduced precision \citep{gra}. Hence, we could implement the matrix-vector multiplication and the orthonormalization in, for example, half precision using fp16 or bfloat16. Second, we can replace a single GMRES solve by an iterative refinement solve by performing reduced-precision GMRES as an inner solver for the corrections \citep{tur}. Third, we could have started our Newton's method with a low precision, e.g., half precision, and only upgrade precision once Newton's method converges to the current precision. We, however, choose not to pursue these enhancements in this work for a few reasons. First, using a reduced precision for the entire GMRES solve is the easiest to implement but gain the most. Also, our primary goal is to demonstrate that Newton-GMRES method can be done in mixed precision, which does not seem previously to have appeared in the literature. The possible enhancements listed above are out of the scope of the current investigation. Second, though the limiting accuracy and residual would stay the same, the effect of reduced precision on the convergence rate is not fully understood for some of these enhancements. The quadratic convergence may be at stake. Finally, hardware support for half precision is not as widely available in CPUs as in GPUs at the time of writing. Moreover, \textsc{Julia} currently supports half precision only through software emulation at, inevitably, a huge cost of speed. We save this line of research for the future.

\subsection{Krylov subspace acceleration}
GMRES also stands a chance for further acceleration with various Krylov 
subspace techniques. In \citet{par}, the authors suggest that the cost of 
constructing the Krylov subspace in an iteration of the Newton-GMRES method 
can be reduced by salvaging the Krylov subspaces of the previous iterations. 
However, our experiments show that the gain acquired from this Krylov subspace 
recycling strategy is very marginal in the current context, if at all, as the Krylov subspaces vary very quickly across Newton iterations. 

Recently, fast randomized sketching is utilized to speed up the subspace projection in GMRES \citep{nak}. However, we found that the sketched GMRES (sGMRES) can hardly accelerate the INGU method. As pointed out in \citet[Section 8.2]{nak}, sGMRES can hardly give a speed boost if the matrix-vector multiplication outweighs other parts of the algorithm, which is exactly the case in the INGU method --- though it is done by FFTs, the Jacobian-vector multiplication is still the most costly part. Therefore, sGMRES has a very large break-even in our experiments and no gain in speed is seen if $n \lesssim 10^4$. Unless a problem is singularly perturbed (see \cref{sec:singular}), there is usually no need for such large degrees of freedom in a 1D problem. 

Meanwhile, sGMRES is unable to handle ill-conditioned systems even when GMRES succeeds \citep[Section 8.2]{nak}. This also rules out the possibility of applying it to singularly perturbed problems, since the condition number of such a problem is usually extremely large. 

Hence, these Krylov subspace acceleration techniques are not included in the experiments shown in \cref{sec:experiments}. However, we believe that these techniques have much bigger potentials for problems in higher spatial dimensions, particularly for those not singularly perturbed.


\section{Numerical experiments} \label{sec:experiments}

\begin{table}[htp]\footnotesize
\begin{minipage}{\linewidth}
\caption{A collection of 1D nonlinear ODE boundary value problems} \label{tab:collection}
\centering
\begin{tabular}{lll}
  \toprule
  equation and BCs & linearization & note \\
  \midrule
  Blasius equation & & \multirow{3}{2.8cm}{boundary layer} \\
  $u''' + uu''/2 = 0 $ & $ \delta''' + (u\delta'' + u''\delta)/2 = 0 $ &  \\
  $ u(0) = 0, u'(0) = 0, u'(L) - 1 = 0 $ & $ \delta(0) = 0, \delta'(0) = 0, \delta'(L) = 0 $ & \\
  \hline
  \addlinespace[1pt]
  Falkner-Skan equation & & \multirow{3}{2.8cm}{an extension of the Blasius equation} \\
  $u''' + uu''/2 + 2\left( 1-(u')^2 \right)/3 = 0 $ & $ \delta''' + (u\delta'' + u''\delta)/2 - 4 u'\delta'/3 = 0$ &  \\
  $ u(0) = 0, u'(0) = 0, u'(L) - 1 = 0 $ & $ \delta(0) = 0, \delta'(0) = 0, \delta'(L) = 0 $ & \\
  \hline
  \addlinespace[1pt]
  Fisher-KPP equation & & \multirow{3}{2.8cm}{a perturbed reaction-diffusion equation} \\
  $u'' + u(1-u) = 0 $ & $ \delta'' + \delta - 2u\delta = 0 $ & \\
  $ u(-4) - 1 = 0, u(4) = 0 $ & $  \delta(-4) = 0, \delta(4) = 0 $ & \\
  \hline
  \addlinespace[1pt]
  fourth-order equation & & \multirow{4}{2.8cm}{the equation of highest-order in this collection} \\
  $u^{(4)} - u'u'' + uu''' = 0$ & $ \delta^{(4)} - u'\delta'' - u''\delta' + u\delta''' + u'''\delta = 0 $ &  \\
  $ u(0) = 0,u'(0) = 0,$ & $ \delta(0) = 0,\delta'(0) = 0,$ & \\
  $ u(1) - 1 = 0,u'(1) + 5 = 0 $ & $ \delta(1) = 0,\delta'(1) = 0 $ & \\
  \hline
  \addlinespace[1pt]
  Bratu equation & & \multirow{3}{2.8cm}{no solution when $\beta > 0.878$ \& closed-form solution exists} \\
  $u'' + \beta e^u  = 0 $ & $ \delta'' + \beta e^u \delta = 0 $ & \\
  $ u(-1) = 0, u(1) = 0 $ & $  \delta(-1) = 0, \delta(1) = 0 $ & \\
  \hline
  \addlinespace[1pt]
  Lane-Emden equation & & \multirow{3}{2.8cm}{an IVP solved as a BVP and closed-form solution exists} \\
  $xu'' +2u' + xu^5 = 0 $ & $x\delta'' +2\delta' + 5xu^4\delta = 0$ &  \\
  $ u(0) - 1 = 0, u'(0) = 0 $ & $  \delta(0) = 0, \delta'(0) = 0 $ & \\
  \hline
  \addlinespace[1pt]
  gulf stream & & \multirow{3}{2.8cm}{a conservation law holds for $u$}\\ 
  $u''' - \beta \left(\left(u'\right)^2-uu''\right) - u + 1 = 0$ & $\delta''' - \beta \left(2u'\delta'-u\delta'' - u''\delta \right) - \delta = 0 $ & \\
  $u(0) - 1 = 0, u'(0)=0, u(L) - 1 = 0$ & $ \delta(0) = 0, \delta'(0)=0, \delta(L) = 0  $ & \\
  \hline
  \addlinespace[1pt]
  interior layer & & \multirow{3}{2.8cm}{singularly perturbed by the leading coefficient} \\
  $\epsilon u'' + uu' + u = 0 $ & $ \epsilon \delta'' + u\delta' + u'\delta + \delta = 0 $ & \\
  $ u(0) + 7/6 = 0 , u(1) - 3/2 = 0 $ & $  \delta(0) = 0, \delta(1) = 0 $ & \\
  \hline
  \addlinespace[1pt]
  boundary layer & & \multirow{3}{2.8cm}{ditto} \\
  $\epsilon u'' + uu' - xu = 0 $ & $ \epsilon \delta'' + u\delta' + u'\delta -x\delta = 0 $ & \\
  $ u(0) + 7/6 = 0 , u'(1) - 3/2 = 0 $ & $  \delta(0) = 0, \delta'(1) = 0 $ & \\
    \hline
  \addlinespace[1pt]
  sawtooth & & \multirow{3}{2.8cm}{ditto} \\
  $\epsilon u'' + (u')^2 - 1 = 0 $ & $\epsilon \delta'' + 2u'\delta' = 0 $ & \\
  $ u(-1) - 0.8 = 0, u(1) - 1.2 = 0 $ & $  \delta(-1) = 0, \delta(1) = 0 $ & \\
  \hline
  \addlinespace[1pt]
  Allen-Cahn equation & & \multirow{3}{2.8cm}{singularly perturbed steady state equation} \\
  $\epsilon u'' + u - u^3 - \sin(x) = 0 $ & $\epsilon \delta'' +\delta - 3u^2 \delta = 0$ & \\
  $u(0) - 1 = 0, u(10) + 1 = 0$ & $  \delta(0) = 0, \delta(10) = 0  $ & \\
  \hline
  \addlinespace[1pt]
  pendulum & & \multirow{3}{2.8cm}{multiple solutions} \\
  $u'' + \sin u  = 0 $ & $ \delta'' + \cos u \delta = 0 $ & \\
  $ u(0) - 2 = 0, u(10) - 2 = 0 $ & $  \delta(0) = 0, \delta(10) = 0 $ & \\
  \hline
  \addlinespace[1pt]
  Carrier equation & & \multirow{3}{2.8cm}{ditto} \\
  $\epsilon u'' +2(1-x^2)u + u^2 - 1 = 0 $ & $\epsilon \delta'' +2(1-x^2)\delta + 2u\delta = 0$ &  \\
  $ u(-1) = 0, u(1) = 0 $ & $  \delta(-1) = 0, \delta(1) = 0 $ & \\
  \hline
  \addlinespace[1pt]
  Painlev\'e equation & & \multirow{3}{2.8cm}{ditto} \\
  $u'' - u^2 + x = 0 $ & $\delta'' - 2u\delta = 0$ & \\
  $ u(0) = 0, u(L) - \sqrt{L} = 0$ & $  \delta(0) = 0, \delta(L) = 0  $ & \\
  \hline
  \addlinespace[1pt]
  Birkisson I & & \multirow{3}{2.8cm}{closed-form solution exists} \\
  $u'' - (\cos x)u' + u\log u = 0$ & $\delta'' - (\cos x)\delta' + (\log u + 1)\delta = 0$ &  \\
  $ u(0) - 1 = 0, u\left(\pi / 2 \right) - e = 0 $ & $  \delta(0) = 0, \delta\left(\pi / 2 \right) = 0  $ & \\
  \hline
  \addlinespace[1pt]
  Birkisson II & & \multirow{3}{2.8cm}{ditto} \\
  $u'' - u' + e^{2x}u + u^2 = \sin^2\left(e^x\right)$ & $\delta'' - \delta' + e^{2x}\delta + 2u\delta = 0$ & \\
  $ u(0) - \sin 1 = 0, u\left(5/2\right) - \sin \left(e^{5/2}\right) = 0$ & $  \delta(0) = 0, \delta\left(5 / 2\right) = 0 $ & \\
  \hline
 \addlinespace[1pt]
  Birkisson III & & \multirow{3}{2.8cm}{ditto} \\
  $u'' + 18(u-u^3) = 0$ & $\delta'' + 18(\delta-3u^2\delta) = 0$ & \\
  $ u(-1) + \tanh 3 = 0, u(1) - \tanh 3 = 0 $ & $  \delta(-1) = 0, \delta(1) = 0  $ & \\
  \bottomrule
\end{tabular}
\end{minipage}
\end{table}

We collect $17$ univariate nonlinear ODE boundary value problems from various sources, such as \citet{bir2, dri}, and use them as a test bank for the INGU framework. These problems and their linearizations are gathered in \cref{tab:collection} with a one-liner note for each problem. For the singularly-perturbed problems, the parameter $\epsilon$ is set to their default values as in \citet{bir2, dri} (see also \cref{tab:data}), only except in \cref{sec:singular} where we focus on solving these equations with small $\epsilon$. Five of these problems have a closed-form solution: the solutions to the Bratu and Lane-Emden equations can be found in standard references of applied mathematics, while those to the Birkisson equations are given in \citet{bir2}. The residual and the error are measured in the $2$-norm and we choose to have GMRES restarted after $r = n/100$ iterations. In case that $r<20$ or $r>150$, we simply take $r = 20$ and $r = 150$, respectively. Except the experiments of \texttt{chebop} which is done with \textsc{Matlab} R2021a, all the numerical experiments are performed in \textsc{Julia} v1.8.2 on a desktop with a 6 core 4.0 Ghz Intel Core i5 CPU and 16GB RAM.


\subsection{Preconditioning} \label{sec:eigen}
\begin{figure}[!ht]
\centering
\subfloat[Blasius equation]{
\includegraphics[width = 11.15cm]{"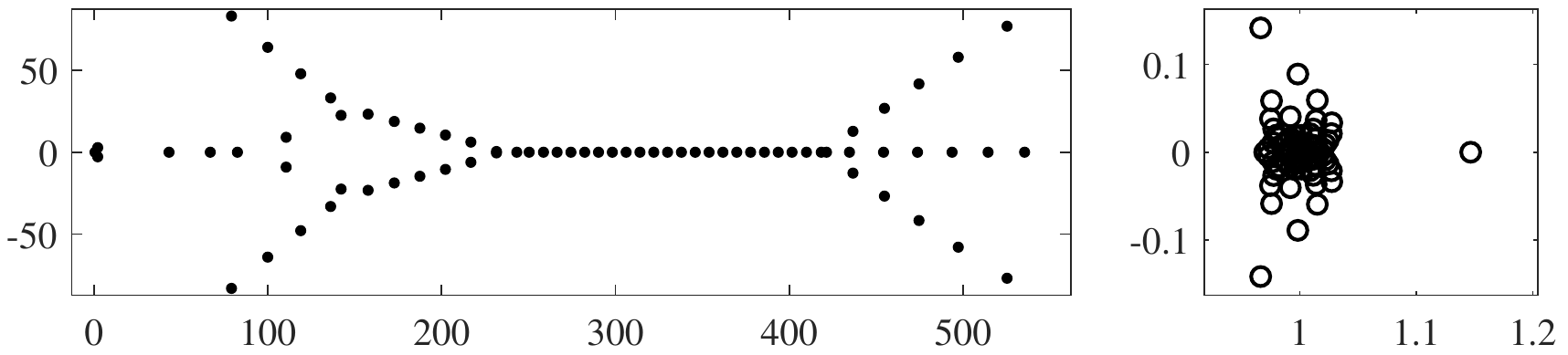"}}
\quad
\subfloat[fourth-order equation]{
\includegraphics[width = 11cm]{"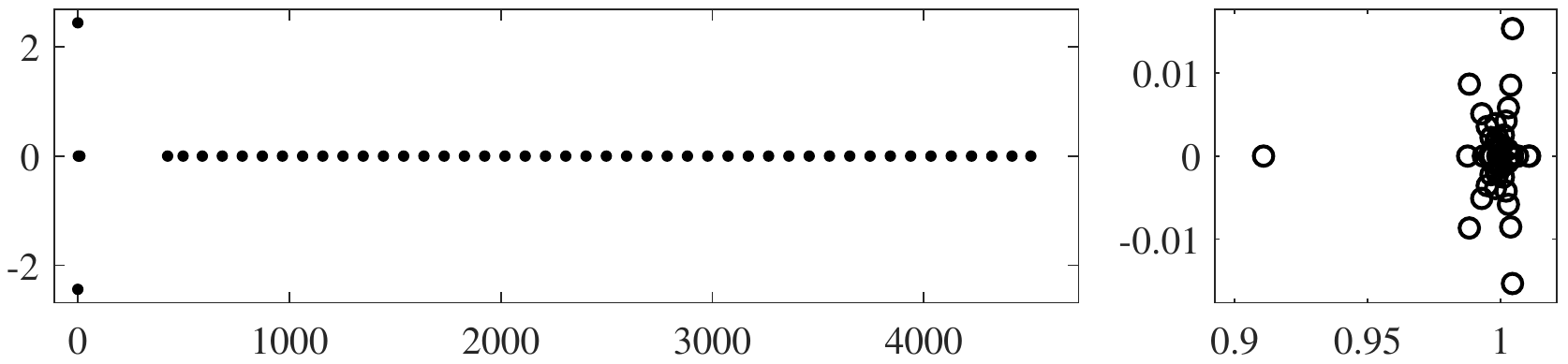"}}
\quad
\subfloat[Lane-Emden equation]{
\includegraphics[width = 11cm]{"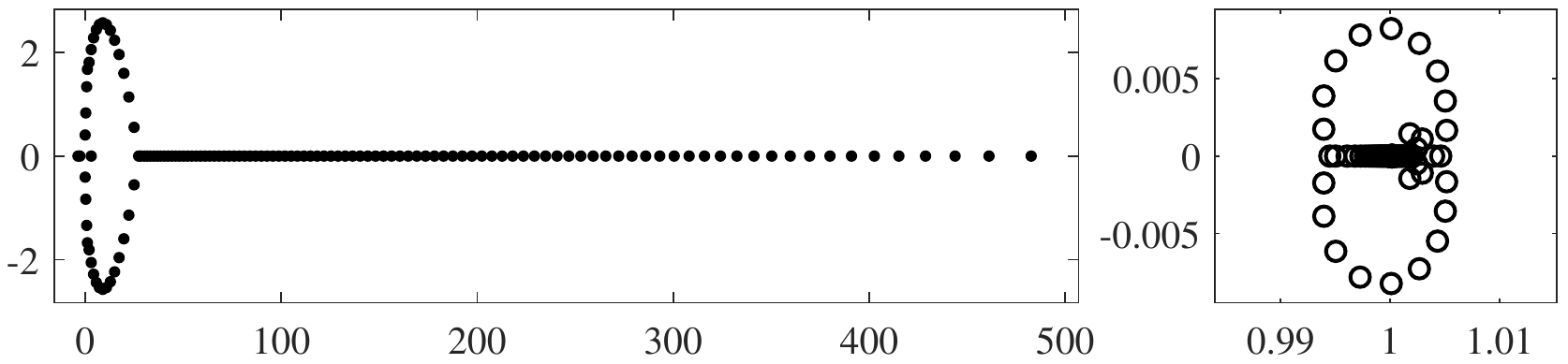"}}
\quad
\subfloat[interior layer equation ($\epsilon = 0.01$)]{
\includegraphics[width = 11cm]{"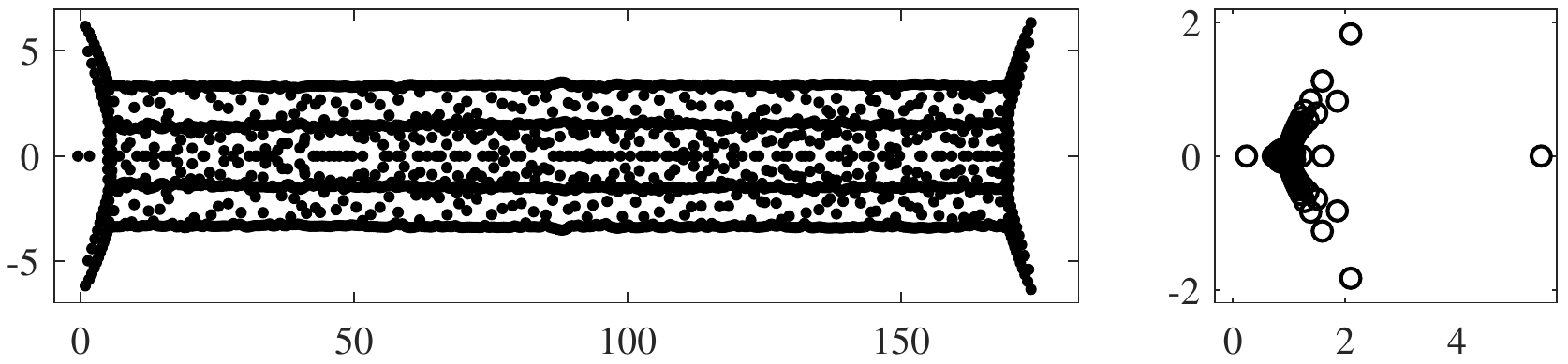"}}
\caption{The eigenvalue distribution of the original systems (left panes) and the preconditioned ones (right panes) for the four selected problems. Note the marked contrast of the axis scales between the left and the right panes of each pair.}
\label{fig:eigen}
\end{figure}
To demonstrate the effectiveness of the preconditioner, we solve the Blasius, 
the forth-order, the Lane-Emden and the interior layer equations and compare 
the eigenvalue distributions of $J_n^k$ with those of $J_n^k 
\left(W_n^k\right)^{-1}$ side by side in \cref{fig:eigen} for the systems arise 
in the last intermediate iteration of each problem. These four problems are 
typical and representative of the entire collection. The Blasius equation 
features a mild boundary layer formed by a physically meaningful no-slip 
boundary condition rather than a boundary layer caused by singular 
perturbation. The fourth-order equation is the one of highest order in this 
collection. The Lane-Emden equation is an ODE initial value problem which we 
solve by regarding it as an ODE boundary value problem. A closed-form 
solution of the Lane-Emden equation is known, which helps in measuring the 
solution error instead of the residual. The equation labeled as interior layer 
has an interior layer caused by singular perturbation in the leading order term 
with $\epsilon = 0.01$ and this interior layer results in a solution with 
length greater than $1000$ for a complete resolution.

The stark contrast between the axis scales of the left panes for the original systems and those of the right panes for the preconditioned shows the extent to which the conditioning is improved. The eigenvalues of the preconditioned system cluster about unity in the complex plane. 

\subsection{Convergence} \label{sec:convergence}
\begin{figure}[t!]
\centering
\subfloat[Blasius equation]{
    \includegraphics[width = 5.5cm]{"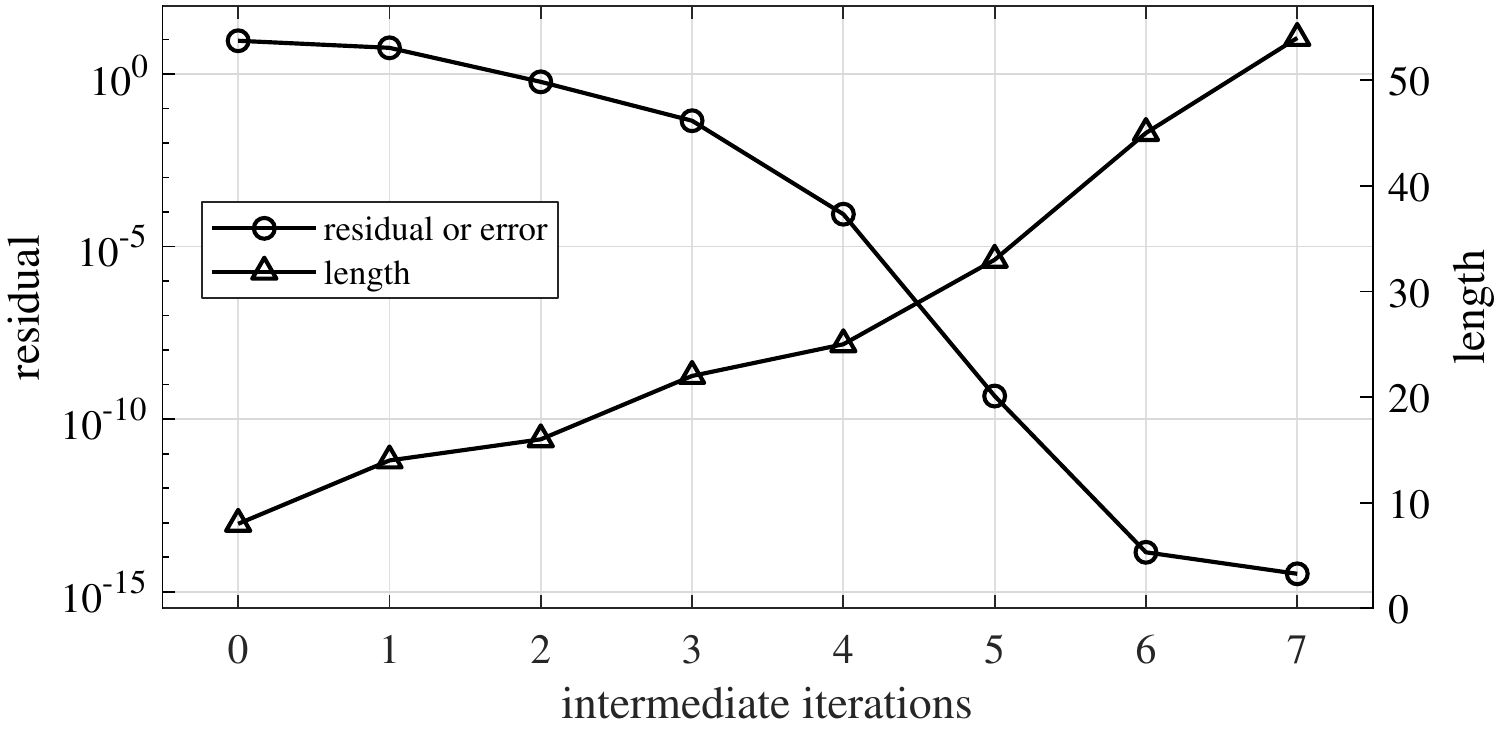"}
  }
  \quad
  \subfloat[fourth-order equation]{
    \includegraphics[width = 5.5cm]{"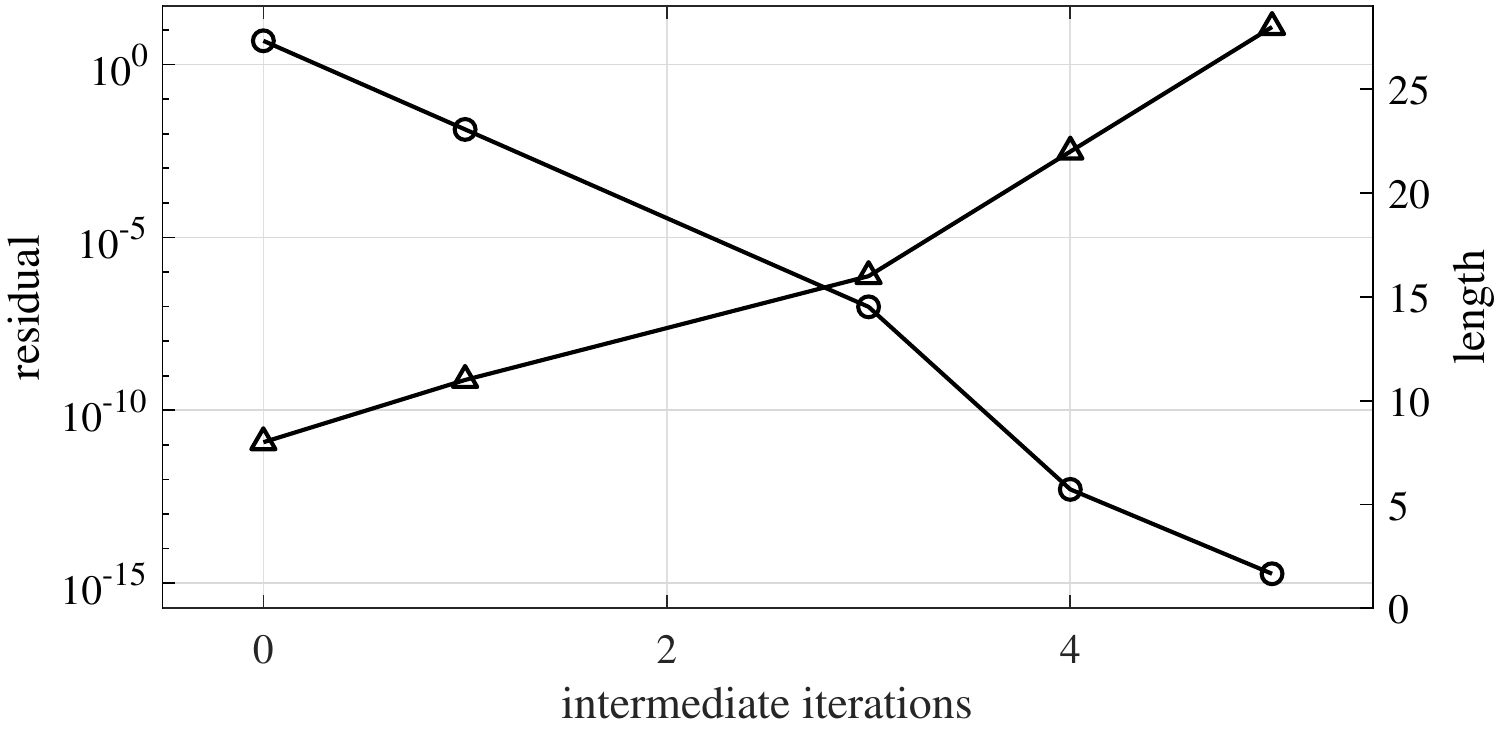"}
  }
  \quad
  \subfloat[Lane-Emden equation]{
    \includegraphics[width = 5.5cm]{"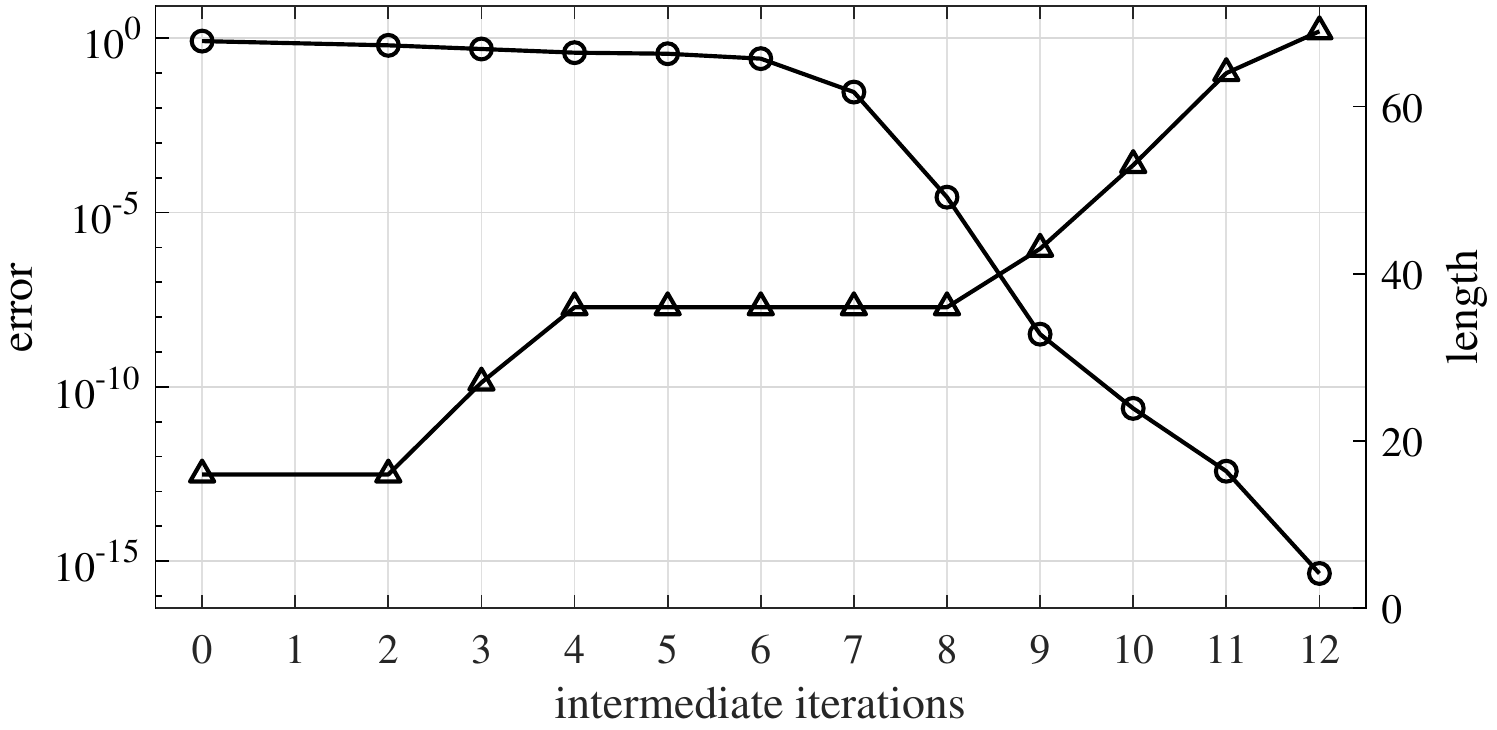"}
  }
  \quad
  \subfloat[interior layer equation ($\epsilon = 0.01$)]{
    \includegraphics[width = 5.5cm]{"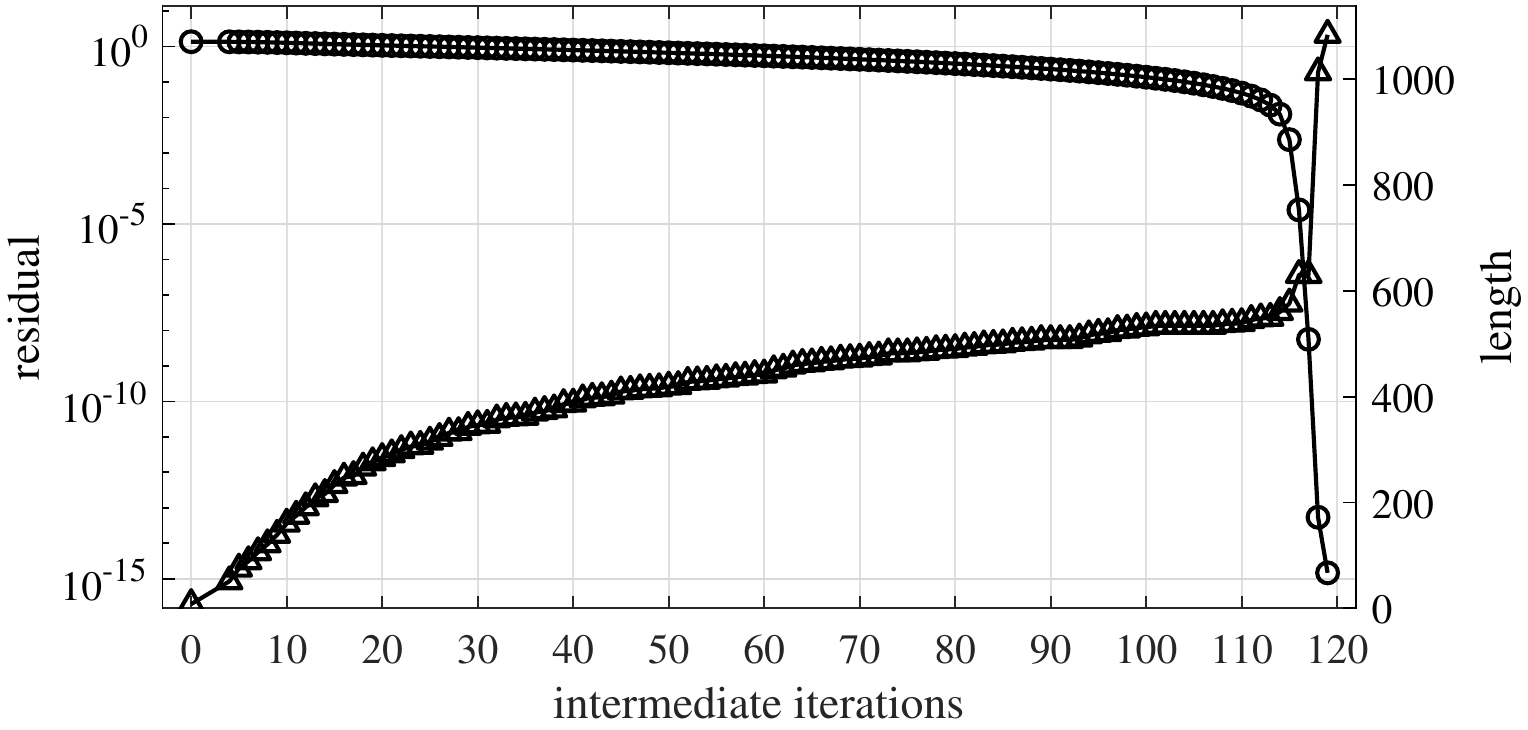"}
  }
\caption{Convergence (vertical axis on the left) and the increase of the solution length (vertical axis on the right) for the four selected problems.}
\label{fig:convergence}
\end{figure}
We take the four problems above again for examples to show the typical 
convergence of the INGU framework. In \cref{fig:convergence}, we show the 
convergence history by plotting the residual $\mG(u^k)$ or absolute error and 
the lengths of the approximate solutions for the TR-contravariant method. The 
other two global methods produce very similar results. The values of the 
residual/error and the solution length can be read off from the $y$-axes on the 
left and the right, respectively. The $x$-axis indicates the number of 
intermediate iterations. A marker signals the start of an outer iteration.

The residual curves for all four problems show that the convergence usually 
evolves with two phases -- the first phase is characterized by the relatively 
level trajectory which corresponds to the slow convergence in the global Newton 
stage and the second phase features a steep descent of the residual/error which 
is the consequence of the fast convergence of the local Newton stage. The 
pattern is most obvious for the interior layer problem with the first $114$ 
intermediate iterations for the global Newton search until the quadratic 
convergence kicks in at about the $5$th to last iteration.

The evolution of the solution lengths matches the decay of the residual in that the lengths increase only moderately until the local Newton phase is reached where the lengths grow very quickly for much improved resolution.

\subsection{Speed and accuracy} \label{sec:speed}
We benchmark the performance of the INGU method against \texttt{chebop} \citep{bir2} by implementing it in conjunction with each of the three global Newton methods of \cref{sec:global}. The \texttt{chebop} class, as a part of the \textsc{Chebfun} system \citep{dri}, is a \textsc{Matlab} solver for linear and nonlinear ODE boundary-value problems. It is extremely well designed and is equipped with sophisticated functionalities of automatic differentiation, linearity detection, lazy evaluation, etc. For nonlinear ODE boundary value problems, \texttt{chebop} employs the trust region method in the affine-variant framework \citep[Section 2.1 \& 3.3]{deu} as the global method and enforces the exact Newton condition. As opposed to TR-contravariant, the affine-variant based trust region method assumes that $\mF(u)$ satisfies an affine variant Lipschitz condition based on which the minimization of the residual is modeled as a constrained quadratic optimization problem. Thus, unlike TR-contravariant, \texttt{chebop} is error-oriented instead of residual-oriented. When the `coefficient' mode is selected, \texttt{chebop} solves the linearized equation using the ultraspherical spectral method. The resulting matrices are stored as sparse data sets, and the \textsc{Matlab} backslash is used for the linear solve. However, it is the LU solver for dense matrices that is usually called, due to the high band density. Therefore, no advantage is taken from the structure of the operators for acceleration.

\begin{figure}[!t]
  \centering
  \subfloat[Blasius equation]{
    \includegraphics[width = 5.5cm]{"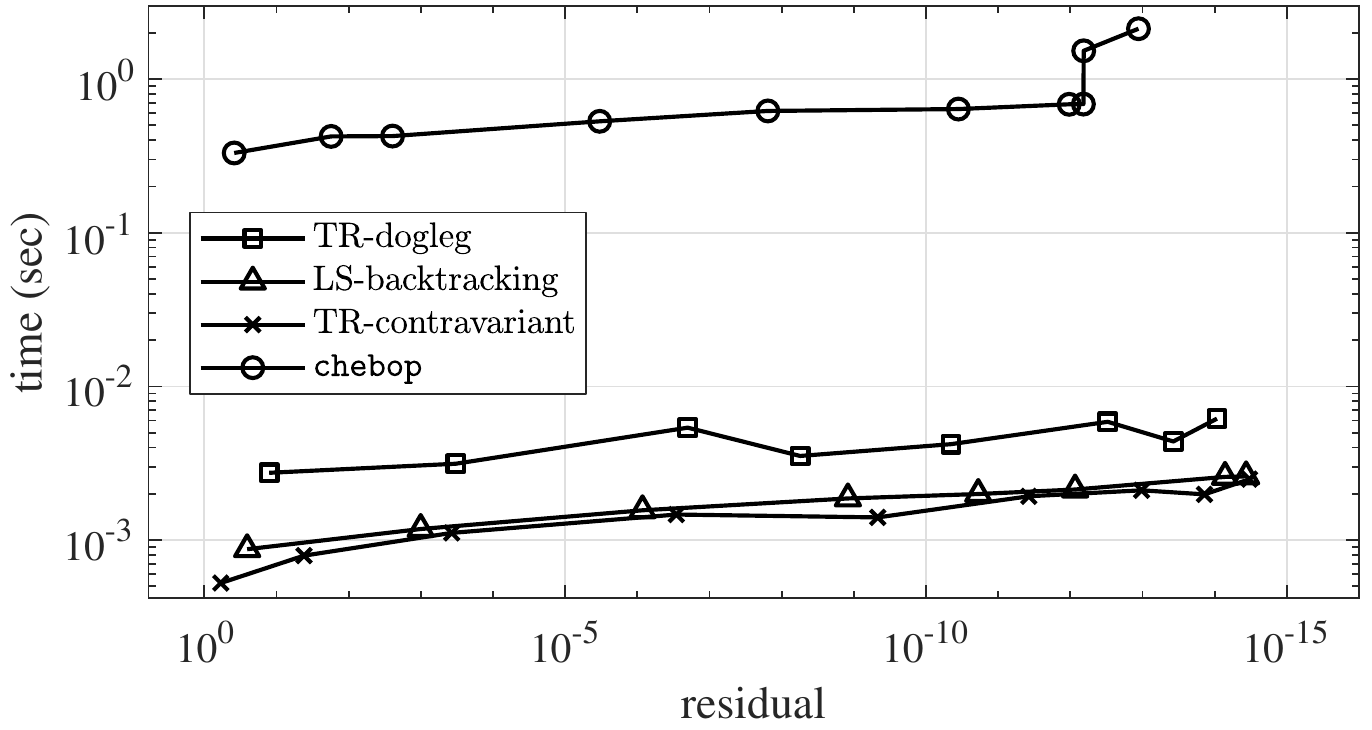"}
  }
  \quad
  \subfloat[fourth-order equation]{
    \includegraphics[width = 5.5cm]{"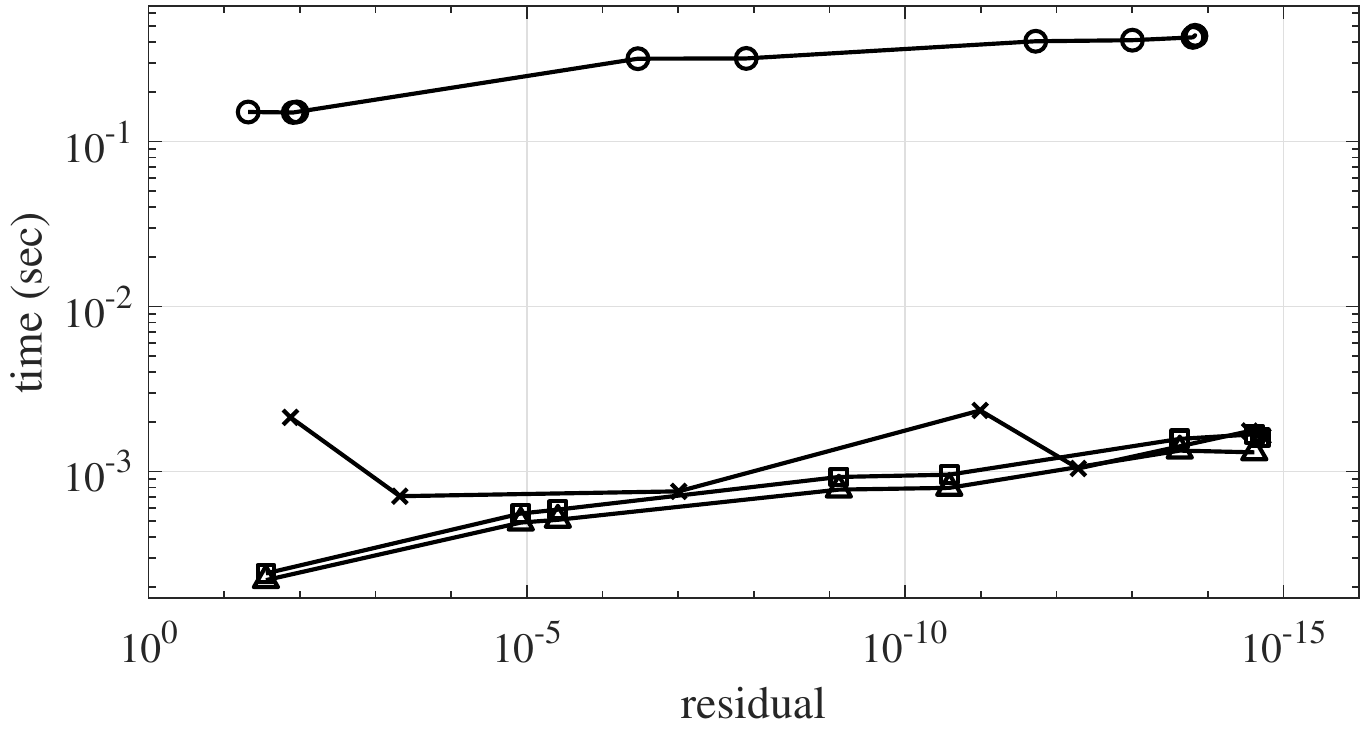"}
  }
  \quad
  \subfloat[Lane-Emden equation]{
    \includegraphics[width = 5.5cm]{"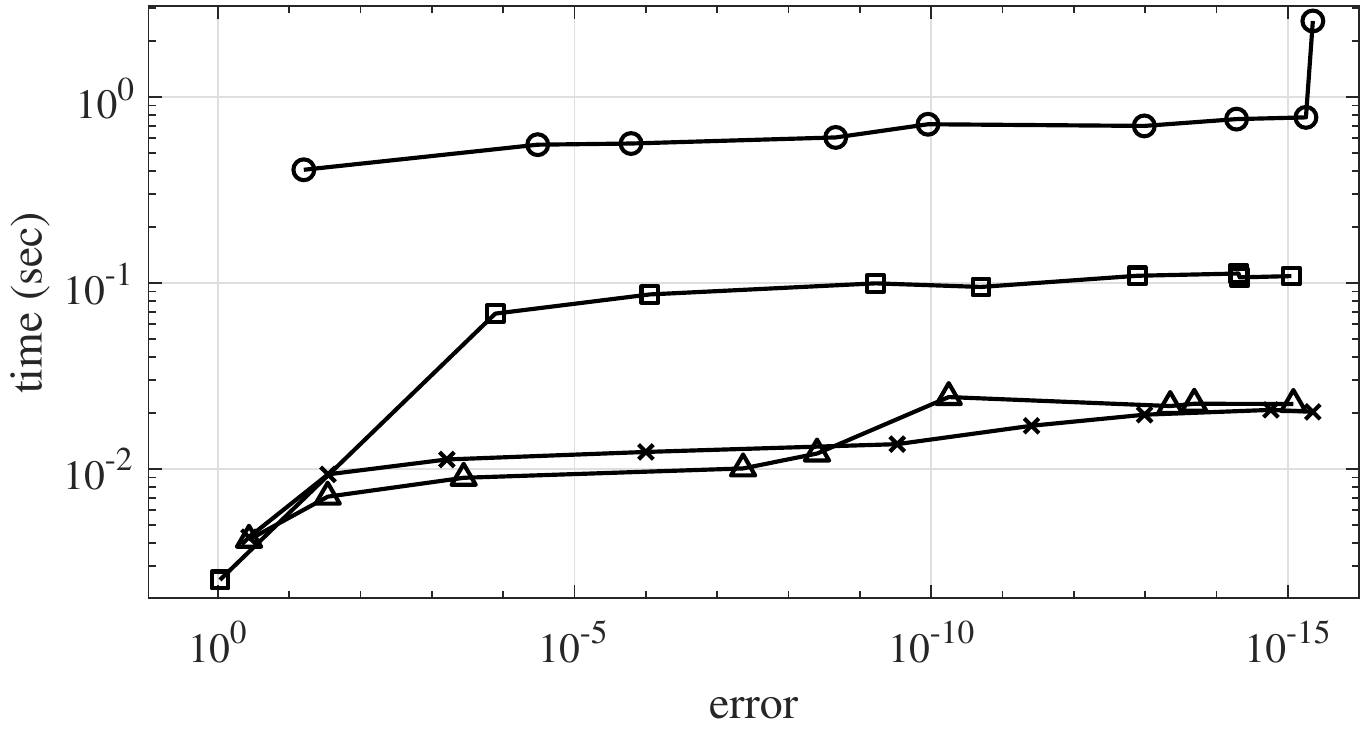"}
  }
  \quad
  \subfloat[interior layer ($\epsilon = 0.01$)]{
    \includegraphics[width = 5.5cm]{"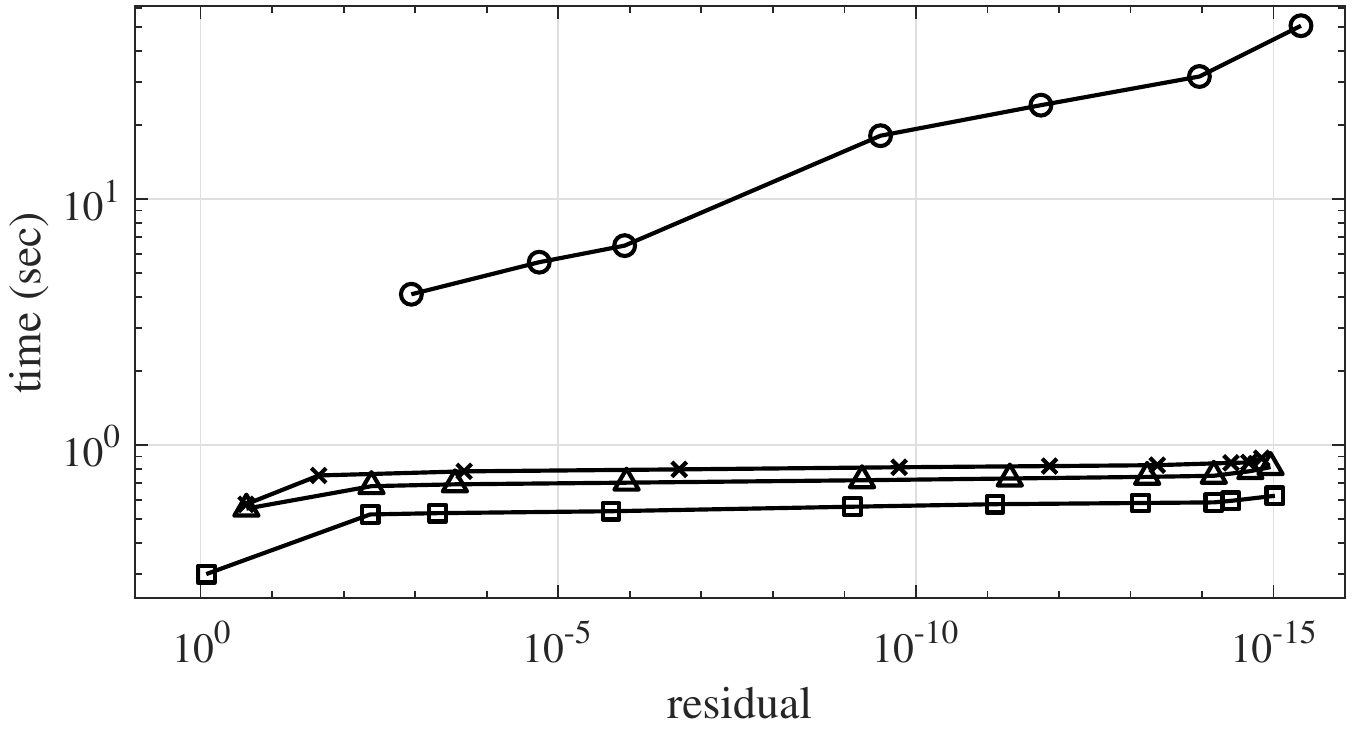"}
  }
  \caption{Residual or error versus execution time for the selected problems.}
  \label{fig:speed}
\end{figure}

\begin{table}[htp]
  \caption{Solution details of the problems in the collection.}\label{tab:data}
  \begin{minipage}{\linewidth}
  \centering
  \begin{tabular}{lcccccc}
  \toprule
  equations &	TRD	&	LSB	&	TRC	&	\texttt{chebop} & speed-up  & length \\
  \midrule
    Blasius & 9.30e-15 & 3.67e-15 & 3.30e-15 & 1.14e-13 & \multirow{2}[0]{*}{8.54e2} & \multirow{2}[0]{*}{54} \\
    ($L=10$) & 6.17e-3 & 2.60e-3 & 2.49e-3 & 2.13 &       &  \\
    \hline
    \addlinespace[1pt]
    Falkner-Skan & 2.91e-14 & 1.65e-14 & 1.98e-14 & 1.77e-13 & \multirow{2}[0]{*}{3.67e2} & \multirow{2}[0]{*}{40} \\
    ($L=10$) & 5.79e-3 & 2.86e-3 & 5.66e-3 & 1.05 &       &  \\
    \hline
    \addlinespace[1pt]
    \multirow{2}[0]{*}{Fisher-KPP} & 1.17e-15 & 1.74e-15 & 1.91e-15 & 7.85e-17 & \multirow{2}[0]{*}{2.33e2} & \multirow{2}[0]{*}{56} \\
          & 4.38e-3 & 3.17e-3 & 2.54e-3 & 5.92e-1 &       &  \\
    \hline
    \addlinespace[1pt]
    \multirow{2}[0]{*}{fourth-order} & 2.03e-15 & 2.42e-15 & 1.86e-15 & 1.43e-14 & \multirow{2}[0]{*}{3.33e2} & \multirow{2}[0]{*}{28} \\
          & 1.60e-3 & 1.31e-3 & 1.64e-3 & 4.36e-1 &       &  \\
    \hline
    \addlinespace[1pt]
    Bratu & 2.44e-15 & 4.66e-15 & 1.22e-15 & 1.78e-15 & \multirow{2}[0]{*}{2.61e2} & \multirow{2}[0]{*}{39} \\
    ($\beta=0.875$) & 4.80e-3 & 3.17e-3 & 3.61e-3 & 8.28e-1 &       &  \\
    \hline
    \addlinespace[1pt]
    \multirow{2}[0]{*}{Lane-Emden} & 8.88e-16 & 8.33e-16 & 4.44e-16 & 4.44e-16 & \multirow{2}[0]{*}{1.26e2} & \multirow{2}[0]{*}{69} \\
          & 1.09e-1 & 2.24e-2 & 2.03e-2 & 2.56 &       &  \\
    \hline
    \addlinespace[1pt]
    gulf stream & 1.43e-12 & 7.62e-13 & 1.43e-13 & 5.28e-13 & \multirow{2}[0]{*}{6.00e1} & \multirow{2}[0]{*}{71} \\
    ($\beta=-0.1,L=35$) & 9.97e-3 & 1.19e-2 & 8.86e-3 & 5.31e-1 &       &  \\
    \hline
    \addlinespace[1pt]
    interior layer & 9.60e-16 & 1.10e-15 & 1.47e-15 & 4.14e-16 & \multirow{2}[0]{*}{8.15e1} & \multirow{2}[0]{*}{1084} \\
    ($\epsilon=0.01$) & 6.22e-1 & 8.14e-1 & 8.85e-1 & 5.07e1 &       &  \\
    \hline
    \addlinespace[1pt]
    boundary layer & 5.62e-15 & 4.34e-15 & 2.75e-15 & 1.83e-12 & \multirow{2}[0]{*}{7.29e1} & \multirow{2}[0]{*}{275} \\
    ($\epsilon=0.01$) & 5.13e-2 & 2.04e-2 & 1.84e-2 & 1.34 &       &  \\
    \hline
    \addlinespace[1pt]
    sawtooth & 6.25e-16 & 3.58e-16 & 3.19e-16 & 6.15e-17 & \multirow{2}[0]{*}{1.15e2} & \multirow{2}[0]{*}{432} \\
    ($\epsilon=0.05$) & 2.15e-2 & 3.06e-2 & 2.48e-2 & 2.49 &       &  \\
    \hline
    \addlinespace[1pt]
    Allen-Cahn & 5.25e-16 & 3.29e-16 & 2.81e-16 & 1.77e-17 & \multirow{2}[0]{*}{8.85e1} & \multirow{2}[0]{*}{79} \\
    ($\epsilon=2$) & 1.37e-2 & 9.95e-3 & 7.30e-3 & 6.46e-1 &       &  \\
    \hline
    \addlinespace[1pt]
    \multirow{2}[0]{*}{pendulum} & 4.59e-15 & 3.66e-15 & 3.70e-15 & 2.06e-15 & \multirow{2}[0]{*}{8.30e1} & \multirow{2}[0]{*}{51} \\
          & 8.92e-3 & 3.33e-3 & 3.94e-3 & 2.76e-1 &       &  \\
    \hline
    \addlinespace[1pt]
    Carrier & 4.60e-16 & 8.52e-16 & 1.89e-16 & 7.86e-17 & \multirow{2}[0]{*}{1.57e2} & \multirow{2}[0]{*}{211} \\
    ($\epsilon=0.01$) & 1.79e-2 & 1.63e-2 & 1.99e-2 & 2.55 &       &  \\
    \hline
    \addlinespace[1pt]
    Painlev\'e & 4.59e-14 & 2.99e-14 & 2.82e-14 & 6.14e-16 & \multirow{2}[0]{*}{9.81e2} & \multirow{2}[0]{*}{51} \\
    ($L=10$) & 5.77e-3 & 2.20e-3 & 2.24e-3 & 2.16 &       &  \\
    \hline
    \addlinespace[1pt]
    \multirow{2}[0]{*}{Birkisson I} & 1.33e-15 & 8.88e-16 & 6.66e-16 & 4.44e-16 & \multirow{2}[0]{*}{9.82e1} & \multirow{2}[0]{*}{23} \\
          & 5.22e-3 & 2.72e-3 & 3.04e-3 & 2.67e-1 &       &  \\
    \hline
    \addlinespace[1pt]
    \multirow{2}[0]{*}{Birkisson II} & 3.33e-15 & 4.11e-15 & 3.22e-15 & 3.22e-15 & \multirow{2}[0]{*}{5.68e1} & \multirow{2}[0]{*}{49} \\
          & 1.43e-2 & 5.84e-3 & 6.37e-3 & 3.32e-1 &       &  \\
    \hline
    \addlinespace[1pt]
    \multirow{2}[0]{*}{Birkisson III} & 9.62e-14 & 1.90e-13 & 4.89e-14 & 1.12e-15 & \multirow{2}[0]{*}{1.46e2} & \multirow{2}[0]{*}{84} \\
          & 1.06e-2 & 6.36e-3 & 1.23e-2 & 9.27e-1 &       &  \\
    \bottomrule
    \end{tabular}   
  \end{minipage}
\end{table}

In \cref{fig:speed}, we plot the accuracy-versus-time curves for the three 
global INGU methods and \texttt{chebop}. The accuracy is measured by the 
absolute residual $\mG(u^k)$ for the approximate solution $u^k$ at the 
termination of the outer iteration or the absolute error in case of an 
closed-form solution being available. We note that the highest accuracy 
achieved by all four methods are much the same. However, it takes significantly 
less time for the INGU methods to achieve a same accuracy goal. Despite the 
different accuracy goals, the speed-up produced by the INGU methods ranges from 
$10 \times$ up to $10^3\times$ compared to \texttt{chebop}. Though our main 
focus is not on the comparison among the three global methods of 
\cref{sec:global}, we find that TR-dogleg is occasionally slower than the other 
two methods. For example, TR-dogleg is about $10\times$ slower than 
TR-contravariant and LS-backtracking for various accuracy goals in the test of 
the Lane-Emden equation. This is partly due to the fact that TR-dogleg also 
involves computation with the transpose of the Jacobians. 

It is worth mentioning that we also reproduced \texttt{chebop} in \textsc{Julia} and find it slightly slower than \texttt{chebop}. This is mainly due to \texttt{SparseArrays}, \textsc{Julia}'s support for sparse vectors, for it is not as efficient as its \textsc{Matlab} counterpart.

Although the observation we make above from the four elaborated examples are representative, we provide solution details in \cref{tab:data} for the entire collection to give a fuller picture of how the INGU framework works. In the first column, we give the values of the parameters, if any, right beneath the equation name. Columns 2-5 contain the minimum residuals or errors that three global INGU implementations can achieve and the corresponding execution time using double precision arithmetic as the working precision. Column 6 lists the speed-up which is the ratio of the \texttt{chebop} execution time to the minimum execution time among the three INGU implementations. Even most of the solutions have a length as small as tens, we gain speed-ups of $100\times$ on average. The rightmost column gives the length of the final solution for the TR-contravariant INGU, as they are typical.

\subsection{Singularly perturbed problems} \label{sec:singular}

So far, we have not yet chosen particularly small $\epsilon$'s for the 
singularly perturbed equations. As $\epsilon$ becomes smaller, difficulty 
arises as the region(s) outside which the initial iterate would fail to converge 
shrinks rapidly. To produce a quality initial iterate, we resort to the 
approach of pseudo-arclength continuation \citep{noc, kel2, bir2}. Specifically, 
we solve the same equation but with a much larger $\epsilon$ for which a simple 
initial iterate usually suffices, e.g., the polynomial of the lowest degree 
that satisfies the boundary conditions. This easy problem can be solved by the 
INGU method up to a low accuracy, and the solution, along with the corresponding 
$\epsilon$, serves as the starting point for the continuation process. It is 
then followed by parameterizing the solution $\displaystyle u(x,s) = 
\sum_{k=0}^{n(s)-1} u_k(s)T_k(x)$ and the perturbation parameter $\epsilon(s)$ 
by the arclength $s$ along the so-called solution path and tracing this path by 
the common predictor-corrector method in order to get close to the target 
$\epsilon$. In each predictor-corrector iteration, we leave the solution path 
by marching along the predictor direction before successive corrector 
iterations bring us (almost) back onto the path so that we have the solution to 
the original equation but with an ever smaller $\epsilon$. The path tracing 
ceases one predictor-corrector iteration before the target $\epsilon$ is overshot. Note that (1) these intermediate solutions usually need not to be 
calculated to high accuracy, as long as they finally lead to a good initial 
iterate; (2) both the predictor and the corrector steps are obtained by solving 
the corresponding expanded equations, and since these expanded equations are 
formed by augmenting \cref{truncated} with one more row and column, the 
solution can be accelerated by the fast multiplication, the preconditioner, and 
the mixed-precision arithmetic as above.

\begin{table}[htp]
  \caption{The intermediate and final solutions to the sawtooth equation on 
the solution path for progressively smaller $\epsilon$.}\label{tab:singular}
    \centering
  \begin{tabular}{crc}
    \toprule
    $\epsilon$ & length & time (sec) \\
    \midrule
    5.00e-2 &		216	 &	1.24e-2\\
    2.52e-2 &		380	 &	8.38e-3\\
    4.30e-3 &		1,299 &		4.77e-2\\
    7.59e-4 &		5,105 &		3.69e-1\\
    1.68e-4 &		21,781 &		2.53\\
    6.97e-5 &		59,741 &		1.45e1\\
    5.33e-5 &		102,748 &		1.95e1\\
    5.05e-5 &		102,748 &		4.09e1\\
    5.00e-5 &		230,755 &		1.40e2\\
    \bottomrule
  \end{tabular}
\end{table}

\begin{figure}[!t]
  \centering
  \subfloat[solutions]{\label{subfig:sol}
    \includegraphics[width = 5cm]{"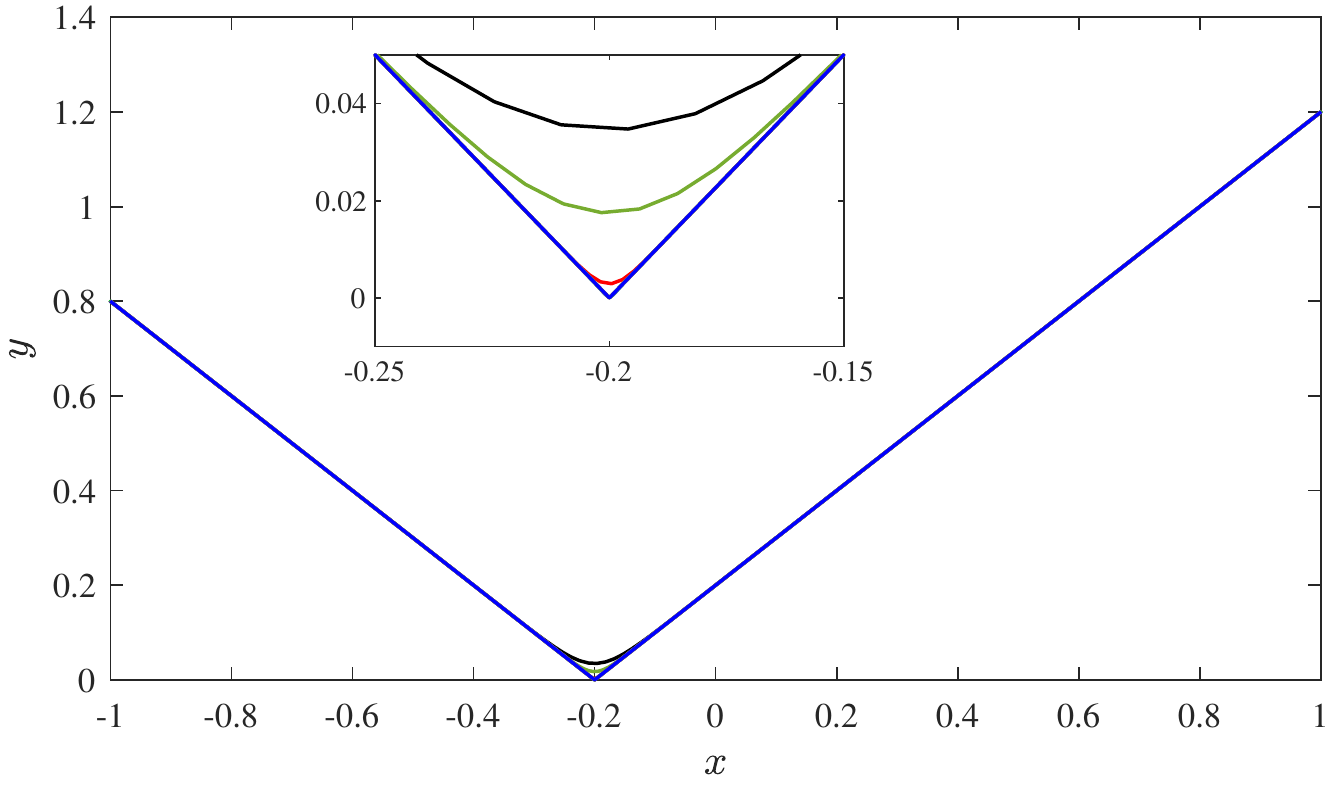"}
  }
  \quad \quad
  \subfloat[coefficients]{\label{subfig:coe}
    \includegraphics[width = 5.5cm]{"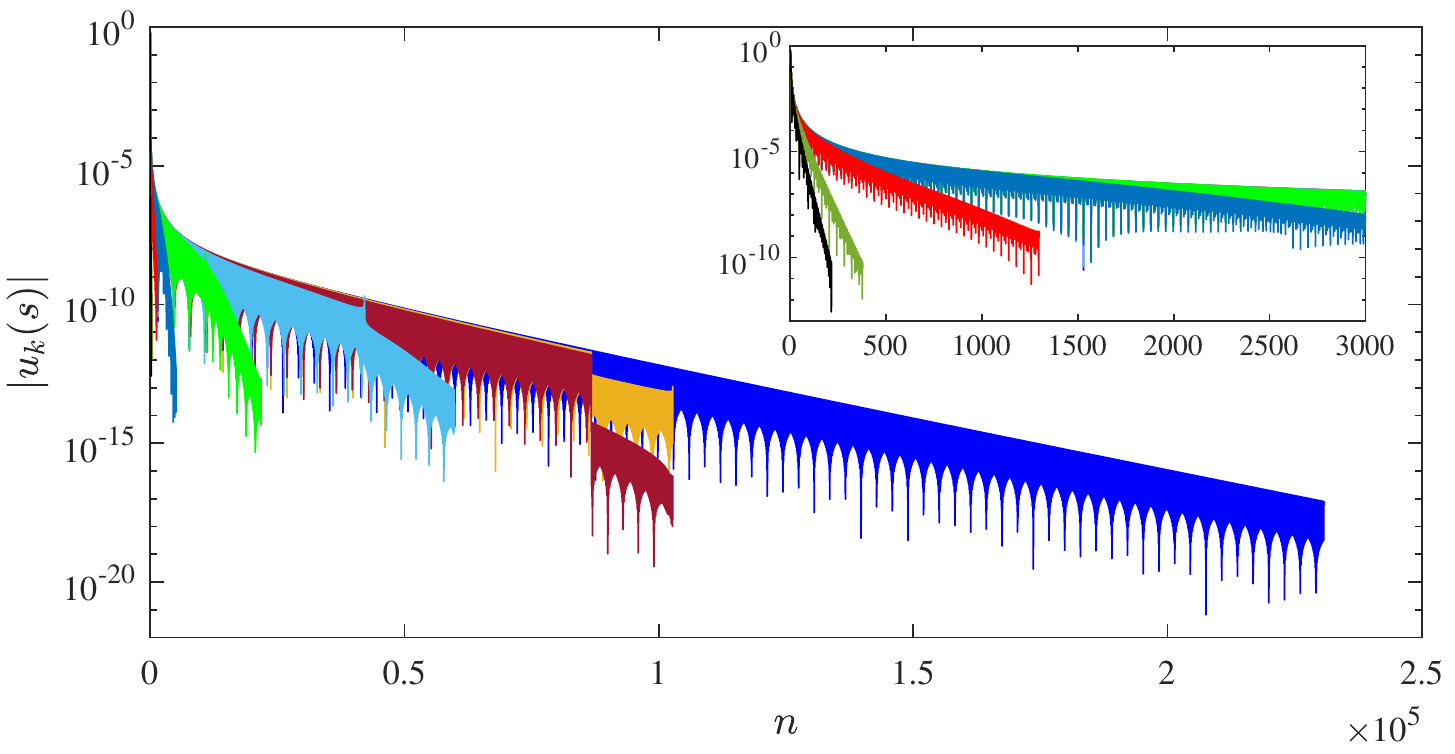"}
  }
  \quad
  \subfloat[execution time and solution length]{\label{subfig:time}
    \includegraphics[width = 5.5cm]{"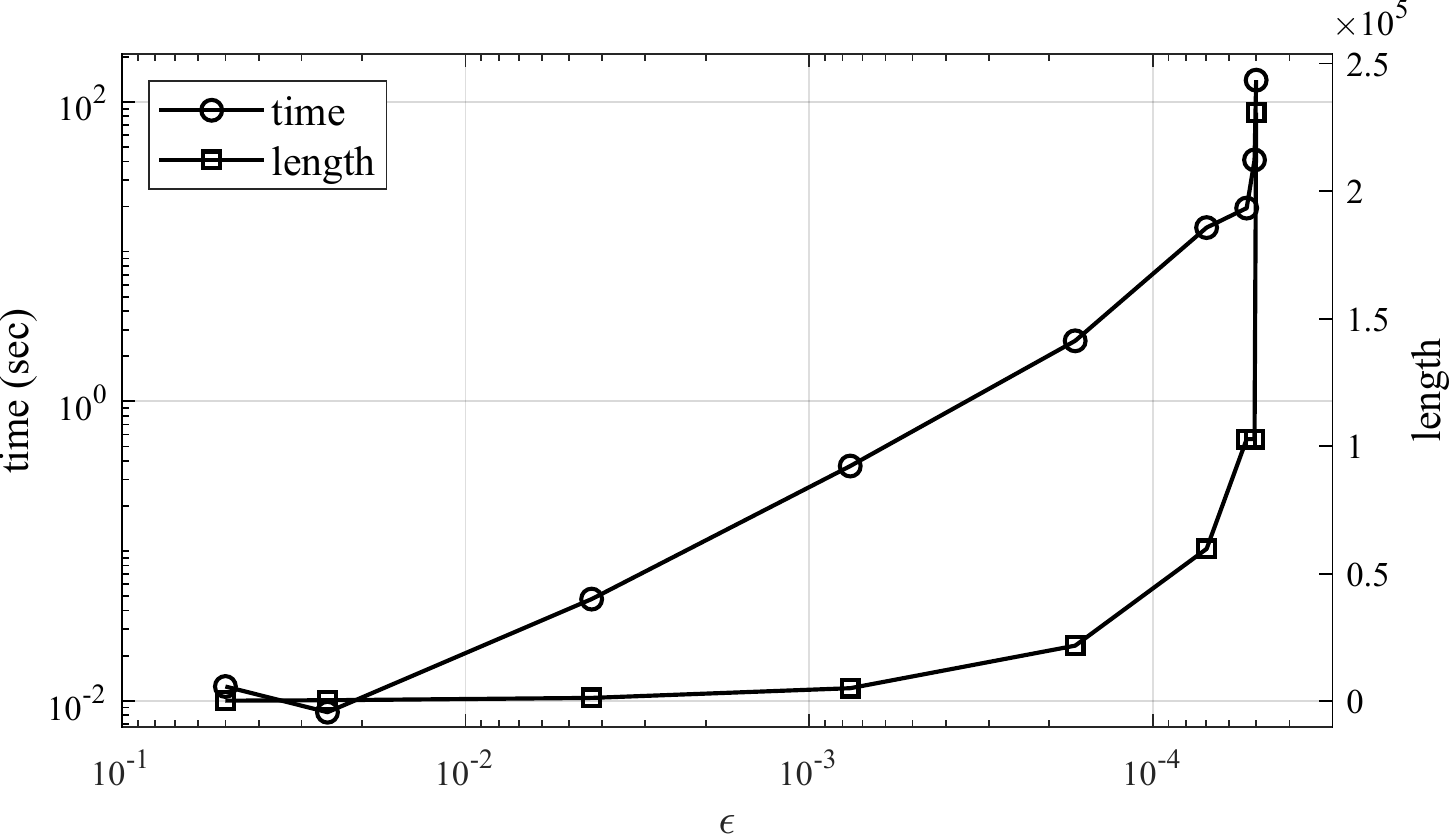"}
  }
  \quad
  \subfloat[solution path]{\label{subfig:path}
    \includegraphics[width = 5.5cm]{"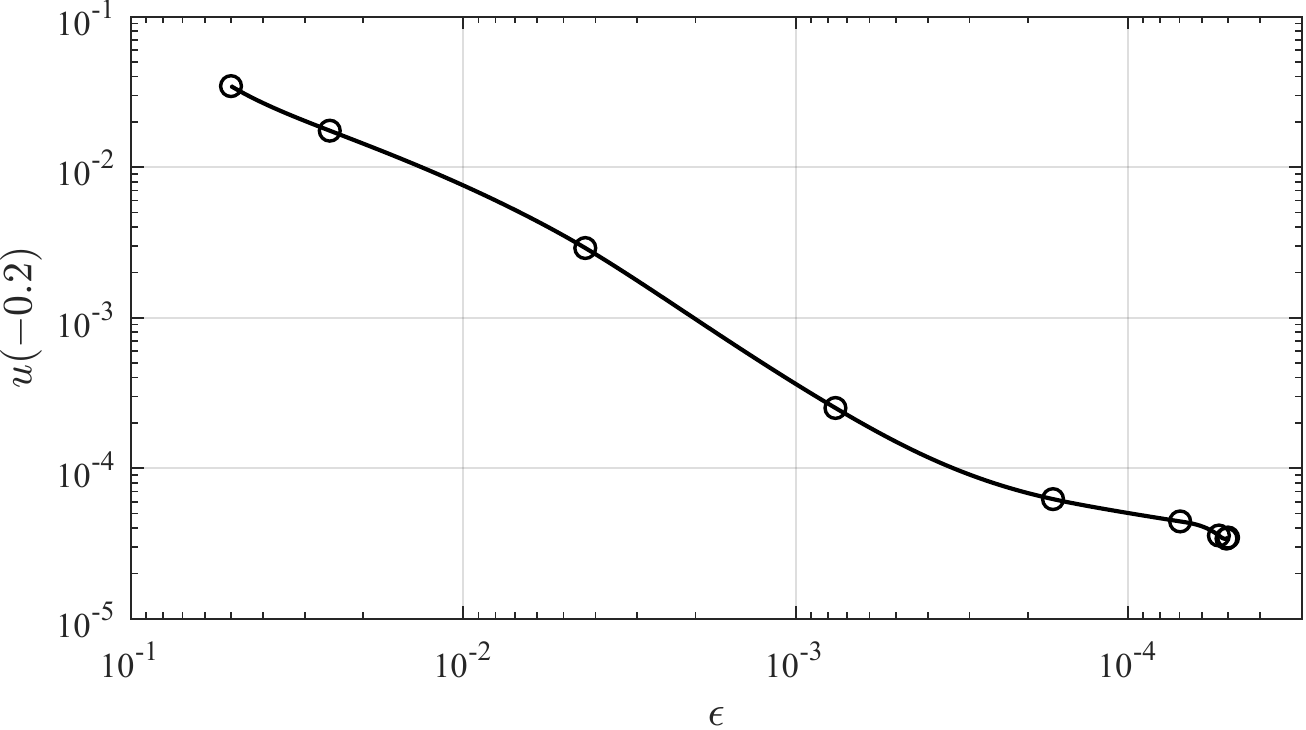"}
  }
  \caption{Solving sawtooth equation with $\epsilon = 5 \times 10^{-5}$.}
  \label{fig:singular}
\end{figure}

We take the sawtooth equation with $\epsilon = 5 \times 10^{-5}$ for example to demonstrate how an average singularly perturbed equation is solved. We start off by solving the sawtooth equation with $\epsilon = 5 \times 10^{-2}$ for which the linear polynomial that satisfies the boundary conditions is good enough as the initial iterate. Setting $s = 0$ for this solution and marching with the predictor-corrector method as described above, we obtain a sequence of 7 more intermediate solutions on the path and the corresponding $\epsilon$'s, before reaching $\epsilon = 5 \times 10^{-5}$. For each predictor step, the corrector iteration stops once the residual of the corrected solution is smaller than $10^{-3}$. See the first 8 rows in \cref{tab:singular} for $\epsilon$ corresponding to these intermediate solutions, their length, and the execution time for computing each of them. These data are also plotted in \cref{fig:singular}\subref{subfig:time}. Finally, the last intermediate solution is fed into the INGU method as the initial iterate. The length of the solution to the target problem and the execution time are appended in the last row of \cref{tab:singular}, and the residual of this final solution is about $1.62\times 10^{-14}$. Note that it takes less than $2.5$ minutes to calculate the final solution whose length is $230,755$! Even if we take into account the pre-computation of the intermediate solutions, the total execution time is still as little as $220$ seconds, which gives a good sense of the speed that the INGU method offers. The scale of the computation in this example almost exhausts the RAM of the machine on which this experiment is carried out. If the memory were large enough, we should be able to solve the sawtooth equation that is much more singularly perturbed.

The sequence of the solutions, including the intermediate ones and the one to the target equation, are shown in \cref{fig:singular}\subref{subfig:sol}, with a close-up displaying the solutions at the turning point $x = -0.2$. As expected, the solutions become more pointed as $\epsilon$ diminishes. The coefficients of these solutions are shown in \cref{fig:singular}\subref{subfig:coe} with a close-up for the first few intermediate solutions that are too short to be seen in the master plot. We also include a plot of the solution path in \cref{fig:singular}\subref{subfig:path} to show the evolution of the value of the solution at $x = -0.2$ versus that of $\epsilon$.

\section{Conclusion} \label{sec:conclusion}
The success of the proposed INGU framework in solving nonlinear equations 
demonstrates that the ultraspherical spectral method is still a powerful tool 
beyond the linear regime, thanks to the structured operators of 
differentiation, conversion, and multiplication.

The INGU framework can readily be extended from ultraspherical or coefficient-based spectral methods for integral equations \citep{sle}, convolution integral equations \citep{xu}, and fractional integral equations \citep{hal} to the nonlinear cases, as the infinite operators in these equations are also structured.

There are both challenges and opportunities towards extending the 
ultraspherical spectral method to solving nonlinear differential equations in 
higher spatial dimensions as generalized Sylvester equations \citep{tow2} and 
tensor decomposition are involved \citep{str}.

\section*{Acknowledgment}
We would like to thank Lu Cheng and Sheehan Olver for their extremely valuable commentary on an early draft of this paper which led us to improve our work. We are also grateful to Joel Tropp for sharing with us his sGMRES code.

\appendix

\section{Post-processing algorithms for the global Newton methods}
We list in this appendix the algorithms of the post-processing step for each of 
the global Newton methods of \cref{sec:global}. Details, including the values 
of the parameters, are given so that the results we report in 
\cref{sec:experiments} can be reproduced when these algorithms are plugged into 
line 7 of \cref{algo:proto}. Particularly, \cref{algo:postTRDogleg} is mainly 
drawn from Algorithm 11.5 and Procedure 11.6 in \citet{noc}. Since no literature 
is found to have a discussion on the determination of the forcing term 
$\omega^{k+1}$ for TR-dogleg, we copy the strategy from the line search method. 
Line 23 and lines 24-25 are borrowed from \citet{eis2} and \citet{kel1}, 
respectively. The main body of \cref{algo:postLSBacktracking} is drawn from 
Algorithm INB in \citet{eis2}, except that lines 21-22 are taken from 
\citet{kel1} and the norm in line 17 is not squared. 
\cref{algo:postTRContravariant} is a combination of the global and local 
inexact Newton-RES methods in \citet[Section 2, 3]{deu}. The description of Newton's 
method used by the \texttt{chebop} nonlinear solver can be found in 
\citet[Algorithm 4]{bir2}. The $2$-norm is used throughout these algorithms.  

\begin{algorithm}[!t]
  \caption{\textsc{Postprocess}: TR-dogleg}
  \label{algo:postTRDogleg}
  \begin{unlist}
      \item[Inputs: ] The current approximate solution $u^k$, nonlinear operator $\mG$, residual vector $f^k$, Jacobian $J_n^k$ and its transpose $\left(J_n^k\right)^{T}$, inexact Newton step $\delta^k$, current size of the trust region $\Delta^k$, termination tolerance $\eta$. 
      \item[Outputs: ] A new approximate solution $u^{k+1}$, a new size of the trust region $\Delta^{k+1}$, a new forcing term $\omega^{k+1}$.
  \end{unlist}
  \hrule
  \begin{algorithmic}[1]
      \State $\bar{\Delta} = 100,~\rho_a = 0.25,~\rho_b = 0.75,~\omega_{\max} = 0.1$  \Comment{initialization}
      \If{$\|\delta^k\| \leq \Delta^k$}
          \State $\widetilde{\delta}^k = \delta^k$  \Comment{inexact Newton step in the trust region}
      \Else
          \State $g^k = \left(J_n^k\right)^{T} f^k$  \Comment{gradient of the merit function}
          \State $\delta^C = - \dfrac{\|g^k\|^2}{\|J_n^k g^k\|^2}g^k$  \Comment{Cauchy point}
          \If{$\|\delta^C\|\geq \Delta^k$}  \Comment{Cauchy point out of the trust region}
              \State $\widetilde{\delta}^k = -\dfrac{\Delta^k}{\|g^k\|} g^k$  \Comment{largest step along Cauchy direction}
          \Else
              \State Find $\nu$ s.t. $\| \delta^C + \nu (\delta^k - \delta^C) \| = \Delta^k$. \Comment{dogleg strategy}
              \State $\widetilde{\delta}^k = \delta^C + \nu (\delta^k - \delta^C)$
          \EndIf
      \EndIf
      \State $\rho = \dfrac{\|f^k\|^2 - \|\mG(u^k + \widetilde{\delta}^k)\|^2}{\|f^k\|^2 - \|f^k + J^k \widetilde{\delta}^k\|^2}$ \Comment{ratio of actual reduction to predicted reduction}
      \If{$\rho < \rho_a$}
          \State $\Delta^{k+1} = \dfrac{\delta^k}{4}$  \Comment{reduction of the trust region}
      \ElsIf{$\rho > \rho_b \And \| \widetilde{\delta}^k \| = \Delta^k$}
          \State $\Delta^{k+1} = \min(\bar{\Delta},~2\Delta^k)$  \Comment{extension of the trust region}
      \Else
          \State $\Delta^{k+1} = \Delta^k$  \Comment{initial size of trust region $\Delta^0 = 0.1$}
      \EndIf
      \If{$\rho > \rho_a$}
      \State $u^{k+1} = u^k + \widetilde{\delta}^k,~f^{k+1} = \mG(u^{k+1})$  \Comment{step accepted}
      \State $\omega^{k+1} = 0.9 \dfrac{\| f^{k+1} \|^{2}}{\| f^k \|^{2}}$  \Comment{initial forcing term $\omega^0=0.1$}
      \State $\omega^{k+1} = \max\left(\omega^{k+1},~\dfrac{\eta}{2\| f^{k+1} \|}\right)$  \Comment{safeguard for oversolving}
      \State $\omega^{k+1} = \min\left( \omega^{k+1},~\omega_{\max} \right)$ \Comment{safeguard for too much inexactness}
      \Else
      \State $u^{k+1} = u^k,~f^{k+1} = f^k,~\omega^{k+1} = \omega^k$  \Comment{no change}
      \EndIf
      \State \Return $u^{k+1}$, $\Delta^{k+1}$, $\omega^{k+1}$
  \end{algorithmic}
\end{algorithm}

\begin{algorithm}[!t]
  \caption{\textsc{Postprocess}: LS-backtracking}
  \label{algo:postLSBacktracking}
  \begin{unlist}
      \item[Inputs: ] The current approximate solution $u^k$, nonlinear operator $\mG$, residual vector $f^k$, inexact Newton step $\delta^k$, current forcing term $\omega^k$, termination tolerance $\eta$. 
      \item[Outputs: ] A new solution $u^{k+1}$ and a new forcing term $\omega^{k+1}$ or a \verb"flag" if fails to converge.
  \end{unlist}
  \hrule
  \begin{algorithmic}[1]
      \State $t = {10}^{-4},~\tau = 1,~\omega = \omega^k,~\kappa=10,~[\gamma_{\min},\gamma_{\max}] = [0.1,0.5],~ \omega_{\max} = 0.1$ \Comment{initialization}
      \State $\omega_s = 0.9 \left(\omega^k\right)^{2}$ \Comment{safeguard}
      \For{$i=1,2,\ldots,\kappa$}  \Comment{finite times of backtracking}
      \State $\tilde{u} = u^k + \tau\delta^k $  \Comment{a trial step}
      \If{$\| \mG(\tilde{u}) \| \leq (1 - t(1-\omega))\| f^k \|$} \Comment{sufficient decrease}
      \State $\widetilde{\delta}^k = \tau\delta^k,~u^{k+1} = u^k + \widetilde{\delta}^k,~f^{k+1} = \mG(u^{k+1})$  \Comment{update}
      \State $\omega^{k+1} = 0.9 \dfrac{\| f^{k+1}\|^{2}}{\| f^k \|^{2}}$  \Comment{initial forcing term $\omega^0 = 0.01$}
      \If{$\omega_s > 0.1$}
      \State $\omega^{k+1} = \max(\omega^{k+1},~\omega_s)$  \Comment{safeguard for too small forcing term}
      \EndIf
      \State $\omega^{k+1} = \max\left(\omega^{k+1},~\dfrac{\eta}{2\| f^{k+1} \|}\right)$  \Comment{safeguard for oversolving}
      \State $\omega^{k+1} = \min\left( \omega^{k+1},~\omega_{\max} \right)$ \Comment{safeguard for too much inexactness}
      \State \Return $u^{k+1}$, $\omega^{k+1}$
      \EndIf
      \State Construct the quadratic polynomial $p(\gamma)$ s.t. $p(0) = g(0),~p'(0) = g'(0),~p(1) = g(1)$, where $g(\gamma) = \| \mG(u^k + \gamma \tau \delta^k) \|$.
      \State Find $\hat{\gamma}$ which minimizes $p(\gamma)$ over interval $[\gamma_{\min}, \gamma_{\max}]$.
      \State $\tau = \hat{\gamma} \tau,~\omega = 1-\hat{\gamma}(1-\omega)$  \Comment{backtracking}
      \EndFor
      \State \Return \verb"flag": failure of backtracking.
  \end{algorithmic}
\end{algorithm}

\begin{algorithm}[!t]
  \caption{\textsc{Postprocess}: TR-contravariant}
  \label{algo:postTRContravariant}
  \begin{unlist}
      \item[Inputs: ] The current approximate solution $u^k$, nonlinear operator $\mG$, residual vector $f^k$, inexact Newton step $\delta^k$, current forcing term $\omega^k$, residual $r^k$ of GMRES, contraction factor $\Theta^{k-1}$ and contravariant Kantorovich quantity $h^{k-1}$ from the previous outer iteration.
      \item[Outputs: ] A new solution $u^{k+1}$ and a new forcing term $\omega^{k+1}$ or a \verb"flag" if fails to converge.
  \end{unlist}
  \hrule
  \begin{algorithmic}[1]
      \State $\mu_{\min}=10^{-6},~\omega_{\max}=0.1,~\omega_{\min}=10^{-5},~\rho = 0.9$  \Comment initialization
      \If{$k \geq 1$}
      \State $\mu = \min\left(1,~\dfrac{1}{(1+\omega^k)\Theta^{k-1} h^{k-1}}\right)$ \Comment{prediction for $\mu$}
      \Else
      \State $\mu = 0.1$
      \EndIf
      \State $\verb"acceptStep" = \verb"FALSE",~\verb"reducted" = \verb"FALSE"$
      \cprotect\While{$ \neg \verb"acceptStep"$}
      \If{$\mu < \mu_{\min}$}
      \State \Return \verb"flag": regularity test fails.
      \Else
      \State $\hat{u}^k = u^k + \mu \delta^k,~\hat{f}^k = \mG(\hat{u}^k)$  \Comment a trial step
      \State $\Theta^{k} = \dfrac{\| \hat{f}^k \|}{\| f^k \|}$ \Comment{contraction factor}
      \State $h^k = \dfrac{2\| \hat{f}^k - (1-\mu)f^k-\mu r^k\|}{\mu^2 \left(1-\left(\omega^k\right)^2\right)\| f^k\|}$ \Comment Kantorovich quantity
      \cprotect\If{$\Theta^k \geq 1 - \frac{\mu}{4}$}   \Comment no contraction
      \State $\mu = \min\left(\dfrac{1}{(1+\omega^k)h^k},~\dfrac{\mu}{2}\right),~\verb"reducted" = \verb"TRUE"$  \Comment damping
      \Else
      \State $\hat{\mu} = \min\left(1,~\dfrac{1}{(1+\omega^k)h^k}\right)$
      \cprotect\If{$\hat{\mu} \geq 4 \mu \And \neg \verb"reducted"$}
      \State $\mu = \hat{\mu}$  \Comment try for a larger step
      \Else
      \State $\verb"acceptStep" = \verb"TRUE"$ \Comment step accepted
      \EndIf
      \EndIf
      \EndIf
      \EndWhile
      \State $\widetilde{\delta}^k = \mu \delta^k,~ u^{k+1} = u^k + \widetilde{\delta}^k,~f^{k+1} = \mG(u^{k+1})$  \Comment update
      \State $\hat{h} = \dfrac{2\rho \left(\Theta^{k}\right)^2}{(1+\rho)\left(1-\left(\omega^k\right)^2\right)}$ \Comment a-posterior estimate
      \State $\omega^{k+1} = \min \left( \dfrac{\sqrt{1+\left(\hat{h}\right)^2} - 1}{\hat{h}},~\omega_{\max} \right)$  \Comment{quadratic convergence mode, $\omega^0=10^{-3}$}
      \State $\omega^{k+1} = \max \left( \omega^{k+1}, \omega_{\min}\right)$  \Comment{safeguard for oversolving}
      \State \Return $u^{k+1}$, $\omega^{k+1}$
  \end{algorithmic}
\end{algorithm}

\bibliographystyle{abbrvnat}
\bibliography{nonlinear-bib}

\end{document}